\documentclass[11pt,twoside]{article}
\usepackage[table,xcdraw]{xcolor} 
\usepackage{amsmath,amsfonts,euscript,amssymb,amsthm,a4,times,color}
\usepackage{microtype}
\usepackage{amssymb,amsmath,amsfonts,amsthm,latexsym,enumerate, esint, manyfoot}
\usepackage{graphics,graphicx}
\usepackage[utf8x]{inputenc}
\usepackage[english]{babel}
\usepackage{amsthm}
\usepackage[notref,notcite]{showkeys}
\usepackage[colorlinks,linkcolor = blue,citecolor = red,pagebackref=false]{hyperref}
\usepackage{amssymb,amsmath,amsfonts,amsthm,latexsym,enumerate, esint, manyfoot}
\usepackage{graphics,graphicx}
\usepackage[utf8x]{inputenc}
\usepackage[english]{babel}
\usepackage{amsmath}
\usepackage{mathabx}
\usepackage{stmaryrd}
\usepackage{mathtools}
\usepackage[notref,notcite]{showkeys}
\usepackage[colorlinks,linkcolor = blue,citecolor = red,pagebackref=false]{hyperref}

\theoremstyle{plain} 
\def\Xint#1{\mathchoice 
  {\XXint\displaystyle\textstyle{#1}}%
  {\XXint\textstyle\scriptstyle{#1}}%
  {\XXint\scriptstyle\scriptscriptstyle{#1}}%
  {\XXint\scriptscriptstyle\scriptscriptstyle{#1}}%
  \!\int} 
\def\XXint#1#2#3{{\setbox0=\hbox{$#1{#2#3}{\int}$} 
  \vcenter{\hbox{$#2#3$}}\kern-.5\wd0}} 
\def\dashint{\Xint-}
\newtheoremstyle{nopunto} 
  {}                       
  {}                       
  {\normalfont}            
  {}                       
  {\bfseries}              
  {}                       
  { }                      
  {\thmname{#1}\thmnumber{ #2}\thmnote{ (#3)}} 
\makeatletter \catcode`@=11
\newbox\tr@tto
\setbox\tr@tto=\hbox{{\count0=0\dimen0=-,9pt\dimen1=1,1pt\loop\ifnum
    \count0<11 \advance \count0 by1 \vrule width.51pt height\dimen1
    depth\dimen0\kern-0.17pt\advance\dimen0 by-0.05pt\advance\dimen1
    by0.1pt\repeat \loop\ifnum\count0<21\advance \count0 by1 \vrule
    width.6pt height\dimen1 depth\dimen0\kern-0.2pt \advance\dimen0
    by-0.1pt\advance\dimen1 by 0.05pt\repeat}}
\def\medint{\displaystyle\copy\tr@tto\kern-10.4pt\int}
\numberwithin{equation}{section}

\allowdisplaybreaks

\theoremstyle{nopunto}
\theoremstyle{plain}
\newtheorem{thm}{Theorem}[section]

\newtheorem{lem}{Lemma}[section]
\newtheorem{prop}{Proposition}[section]
\newtheorem{defn}{Definition}[section]

\usepackage{comment}
\usepackage{xpatch}
\usepackage{orcidlink}
\usepackage{blindtext}

\usepackage{authblk}
\usepackage{hyperref}
\usepackage{url}

\makeatletter
\usepackage{fancyhdr}
\usepackage[export]{adjustbox}
\usepackage{hyperref}
\hypersetup{colorlinks=true, citecolor=cyan, linkcolor=black, urlcolor=black, pdftitle=Tesi, pdfauthor=Anna Cavagnoli}
\usepackage{booktabs}
\usepackage{multirow}

\title{\Large\bfseries Higher differentiability for solutions to Stokes systems under subquadratic growth conditions}

\makeatletter \catcode`@=11
\newbox\tr@tto
\setbox\tr@tto=\hbox{{\count0=0\dimen0=-,9pt\dimen1=1,1pt\loop\ifnum
    \count0<11 \advance \count0 by1 \vrule width.51pt height\dimen1
    depth\dimen0\kern-0.17pt\advance\dimen0 by-0.05pt\advance\dimen1
    by0.1pt\repeat \loop\ifnum\count0<21\advance \count0 by1 \vrule
    width.6pt height\dimen1 depth\dimen0\kern-0.2pt \advance\dimen0
    by-0.1pt\advance\dimen1 by 0.05pt\repeat}}
\def\medint{\displaystyle\copy\tr@tto\kern-10.4pt\int}
\numberwithin{equation}{section}

\allowdisplaybreaks

\begin{document}

\title{Higher differentiability for solutions to stationary $p$-Stokes systems under sub-quadratic growth conditions}

\author{ Anna Cavagnoli \\ 
\normalsize{Dipartimento di Matematica e Applicazioni "R.
Caccioppoli"} \\ \normalsize{Universit\`{a} di Napoli ``Federico
II", via Cintia - 80126 Napoli}
 \\ \normalsize{e-mail:
anna.cavagnoli@unina.it}}

\bigskip

\maketitle

\begin{abstract}
  We consider weak solutions $(u,\pi):\mathbb{R}^n\supset\Omega\to\ \mathbb{R}^n\times\ \mathbb{R}$
  to stationary $p$-Stokes systems of the type
 \[ 
  \begin{cases}
    -\mathrm{div} (a(\mathcal{E} u))+\nabla\pi=f  \\
    \mathrm{div}(u)=0, 
  \end{cases}
  \]
in $\Omega,$ where the function $a(\xi)$ satisfies $p$-growth conditions in
  $\xi$. By $\mathcal{E} u$ we denote
  the symmetric part of the gradient $Du$. In this setting, we establish results on the
  fractional higher
  differentiability of both the symmetric part of the gradient $ D u$
  and of the
  pressure $\pi$.
\end{abstract}
\section{Introduction and statement of the results}
\par
\par
In this paper we study the higher differentiability of solutions $(u,\pi):\Omega\to\ \mathbb{R}^n\times\ \mathbb{R}$
  to stationary $p$-Stokes systems, of the form
\begin{equation} \label{equa}
  \begin{cases}
    -\mathrm{div} (a(\mathcal{E} u))+\nabla\pi = f,  \\
    \mathrm{div}(u) = 0,
  \end{cases}
\end{equation}
in $\Omega,$ a bounded domain in $\mathbb{R}^n$, $ n \geq 2.$ In this context $\mathcal{E}  u $ denotes the symmetric part of the gradient  $ Du \in \mathbb{R}^{n \times n} $ and $f : \Omega \to \mathbb{R}^n $ is a given inhomogeneity. \par
It is worth pointing out that the symmetric gradient of functions plays a key role in the theory of non-Newtonian fluid mechanics, for a detailed description of related models we refer to \cite{FluidMechanics1} and \cite{FluidMechanics2}. In this setting, $u(x)$ and $\nabla u $ represent, respectively, the configuration of a fluid at a point $ x \in \Omega $ and its velocity. The polynomial growth finds a motivation in the non linearity describing the non-Newtonian fluid under consideration. \par
  Here, we investigate the higher differentiability of the solutions $(u,\pi):\Omega\to\mathbb{R}^n\times \mathbb{R}$ under appropriate assumptions on the regularity of the map $\xi\to a(\xi)$ and the integrability of the right-hand side $f$, in case the non-linearity $a(\xi)$ satisfies subquadratic growth conditions. 

As it is customary, due to the nonlinear nature of the problem, it is not to
be expected that second weak derivatives exist, but the extra differentiability can be proven  for a certain nonlinear
quantity of the symmetric gradient, such as $V_\mu(\mathcal{E} u)=(\mu^2 + \left| \mathcal{E}u \right|^2)^\frac{p-2}{4} \mathcal{E} u$. We will prove that this quantity admits a weak derivative, in the sense of  $V_\mu(\mathcal{E} u) \in W^{1,2}_{loc}(\Omega).$
This is consistent with the well-established theory for this special case
of the $p$-Laplace operator, in which higher differentiability results
are known, see \cite{Giusti} and the references therein.

The question of the  higher differentiability  for the
$p$-Stokes system \eqref{equa} has been less investigated. In this case, indeed, more difficulties arise as the system
\eqref{equa} contains only information on the symmetric part of the
gradient.
For Stokes-systems, i.e. without a convective term, corresponding results on
higher differentiability of solutions were established by Naumann
\cite{Naumann} for the case of polynomial growth 
and by Diening \& Kaplick\'y~\cite{DieningKaplicky} under a more
general growth condition of Orlicz type.

All of the mentioned results 
cover only systems with constant coefficients. A more complete picture is contained in \cite{Passarelli:2020:1}, where a
fractional higher differentiability result for $p$-Stokes
and Navier-Stokes systems under 
a H\"older continuity assumption on the  coefficients and with the restrictions
$p>\frac{3n}{n+2}$ and $p\ge2$ has been obtained.\par
The subquadratic growth case has been dealt only as a particular case of the Orlicz growth case and the available results are obtained under assumptions on the datum $f$, that are far from being optimal (See \cite{DieningKaplicky} and \cite{Passarelli:2023:2}). On the other hand, an higher differentiability result for solutions to $p$-harmonic systems with $1<p<2$ has been obtained in \cite{ClopGentPass} under a sharp assumption on the datum $f$ in the scale of Lebesgue spaces. \par
The main goal of this paper is to obtain under the same assumption on $f$ an higher differentiability result of integer order for solutions to Stokes system as in \eqref{equa}.

More precisely, we consider a measurable function $a: \mathbb{R}^{n \times n}_{sym}\to\mathbb{R}^{n \times n}_{sym}$, where
$\mathbb{R}^{n \times n}_{sym}$ denotes the space of symmetric real-valued $n\times n$
matrices,  $n\ge2$.
In what follows, we assume that
$\xi\mapsto a(\xi) \quad \mbox{is $C^1$ }$
and that the following conditions are satisfied for given
parameters $1 < p<2$ and $\mu\in[0,1]$:
\begin{itemize}
    \item there exist positive constants $\ell, L$ such that
\end{itemize}
\begin{equation}\label{ip1}
 \ell(\mu^2+|\xi|^2)^{\frac{p-1}{2}} \leq |a(\xi)| \le L(\mu^2+|\xi|^2)^{\frac{p-1}{2}},
\end{equation}
for a.e. $x \in \Omega $ and every $\xi \in \mathbb{R}^{n \times n}_{sym}.$
\begin{itemize}
    \item there exists a positive constant $\nu > 0 $, such that
    \begin{equation}\label{ip2}
\langle a(\xi) - a(\eta),\xi - \eta \rangle\ge \nu(\mu^2+|\xi|^2 + \left| \eta \right|^2)^{\frac {p-2}{2}}|\xi - \eta|^2,
\end{equation}
\end{itemize}

 for a.e. $ x \in \Omega$ and  every $\xi,\eta\in\mathbb{R}^{n \times n}_{sym}$. 
 \begin{itemize}
     \item there exists a positive constant $L_1 > 0 $ such that
    \begin{equation}\label{ip3}
        \left| a(\xi) - a(\eta) \right| \leq L_1(\mu^2 + \left| \xi \right|^2 + \left| \eta \right|^2)^\frac{p-2}{2} \left| \xi - \eta \right|
    \end{equation}
    for a.e. $x \in \Omega $ and every $ \xi, \eta \in \mathbb{R}^{n \times n}_{sym}.$
 \end{itemize}
  Let us observe that 
  \begin{equation}
      \frac{np}{n(p-1)+2-p} > \frac{np'}{n+p'}= \frac{np}{n(p-1)+p} \notag
  \end{equation}
  if and only if
  \begin{equation}
      n(p-1) +p > n(p-1) + 2-p \notag
  \end{equation}
  but this is true for every $ p > 1. $ For this reason, we can consider a force term 
  \begin{equation}\label{fassum}
      f \in L^\frac{np}{n(p-1)+2-p}_{loc}(\Omega).
  \end{equation}  \par
   Note that $ f \in L^{\frac{np}{n(p-1)+2-p}}_{loc}(\Omega)\subset L^{\frac{np'}{n+p'}}_{loc}(\Omega) \subset W^{-1,p'}_{loc}(\Omega).$\par
We consider weak solutions
$(u,\pi)\in W^{1,p}(\Omega,\mathbb{R}^n)\times L^{p'}(\Omega)$,
in the sense made precise in 
  Definition~\ref{def:weak-solution} below, to the system \eqref{equa}. \par

\begin{defn}\label{def:weak-solution}
We call $(u,\pi)\in W^{1,p}(\Omega,\mathbb{R}^n)\times L^{p'}(\Omega)$ a weak solution of
the system \eqref{equa} if $\mathrm{div} (u)=0$ holds in the sense of
distributions and
\begin{equation*}
  \int_{\Omega}\big(\langle a(\mathcal{E}
  u),\mathcal{E}\varphi\rangle-\pi\,\mathrm{div}\varphi\big)\,dx
  =
  \int_{\Omega}f\cdot\varphi\,dx
  \end{equation*}
  holds for every $\varphi\in C_0^\infty(\Omega,\mathbb{R}^{n})$.\end{defn}

 Since $ f \in W^{-1,p'}_{loc}(\Omega)$, the right-hand side of the previous equation is well-defined.

We use the customary notation
\begin{equation*}
  V_\mu(\xi):= (\mu^2+|\xi|^2)^{\frac{p-2}4}\xi
  \qquad\mbox{for all $\xi\in\mathbb{R}^k$. }
\end{equation*}

Our main result is the following:
\begin{thm}\label{Thm}
    Assume that \eqref{ip1}, \eqref{ip2} and \eqref{fassum} are in force for an exponent $p$ such that $ 1 < p < 2$ and that $ (u, \pi) \in W^{1,p}(\Omega, \mathbb{R}^n) \times L^{p'}(\Omega) $ is a weak solution of the system
    \eqref{equa} in the sense of Definition \ref{def:weak-solution}. Denote with $ q = \frac{np}{n+p-2} $. \par
    Then we have 
    \begin{center}
       $ V_{\mu}(\mathcal{E}u) \in W^{1,2}_{loc}(\Omega, \mathbb{R}^{n \times n}_{sym})  \quad \mathrm{and}  \quad  \pi \in W^{p-1, q'}_{loc}(\Omega).$
    \end{center}
Moreover, the local estimate

    \begin{equation}
\begin{aligned}
\int_{B_{\frac{R}{2}}(x
_0)} \left| D V(\mathcal{E} u) \right|^2 \, dx 
&\leq \frac{c}{R^2} \int_{B_{R}(x_0)} \left( \mu^2 + |\mathcal{E} u|^2 \right)^{\frac{p}{2}} \, dx 
+ \frac{c}{R^2} \int_{B_R(x_0)} \left| D u \right|^p \, dx \\
&\quad + \frac{c}{R^{\frac{2p}{p-1}}} 
\| f \|_{L^{q'}(B_R(x_0
))}
\end{aligned}
\end{equation}
holds true for any ball $ B_R(x_0) \Subset \Omega$ with $ 0 < R < 1$. \par \qquad \\\
For the pressure, we have the following local estimate for every cut-off function $\eta \in C_0^{\infty}(B_{\frac{R}{2}(x_0)}, [0,1])$ such that $\left| \nabla \eta \right| \leq \frac{c}{R} $ 
\begin{align}
  \sup_{0 < \left| h \right| < \frac{R}{4}} \Biggl( \int_{B_R(x_0)} \left| \frac{\tau_h(\eta \pi)}{|h|^{p-1}} \right|^{q'} \, dx \Biggr)^{\frac{1}{q'}} 
    &\leq \|f\|_{L^{q'}(B_R(x_0))} 
    + \left( \int_{B_{R/2}(x
    _0)} |D(\mathcal{E} u)|^p \, dx \right)^{\frac{p-1}{p}} \notag \\
    &\quad + \frac{c}{R} 
    \left( \int_{B_R(x_0)} (\mu^2 + |\mathcal{E} u|^2)^{\frac{np}{2(n-2)}} \, dx \right)^{\frac{(n-2)(p-1)}{np}} 
  \\ + \frac{c}{R} \| \pi \|_{L^{q'}(B_R(x_0))}. 
\end{align}
\end{thm}
\qquad 

Regarding the proof strategy of this Theorem, for the derivation of the higher differentiability result, we apply the difference
quotient method. However, the present situation requires a suitable modification of this technique, because due to the
pressure term in \eqref{equa}, it is only feasible to test \eqref{equa} with divergence-free test functions. This technical difficulty is solved
by means of a well-known lemma by Bogovskiĭ, which is applied to construct a suitable correction term. Once the higher
differentiability of $\mathcal{E}u$ has been established, we can solve the system \eqref{equa} for the pressure and deduce also fractional higher
differentiability of the pressure.\par

Moreover, in this specific case, we work with the symmetric part of the gradient of the solution, which is closely related to the full gradient through a lemma known as the Sobolev–Korn inequality. Some preliminary lemmas, in addition to Lemma \ref{le1}, will be frequently employed to derive upper and lower estimates for $V_{\mu}(\mathcal{E}u)$.
\section{Preliminaries}\label{sec:prelim}
This section is devoted to collect notations and preliminary results that will be needed in what follows.
\subsection{Notation and elementary lemmas}

\noindent

We write $B_\rho(x_0)\subset\mathbb{R}^n$ for the open ball of radius $\rho>0$ and center
$x_0\in\mathbb{R}^n$.
For the mean value of a function $f\in L^1(B_\rho(x_0),\mathbb{R}^k)$, we write
\begin{equation*}
  (f)_{x_0,\rho}:=\dashint_{B_\rho(x_0)}f(x)\,dx
\end{equation*}

For the standard scalar product on the space
$\mathbb{R}^{n\times n}$ of $n\times n$ matrices, we write $\langle\cdot,\cdot\rangle$, in contrast to the
Euclidean scalar product on $\mathbb{R}^n$, which we denote by "$\cdot$".\par
We will  denote by $c$ a
general constant that may vary on different occasions, even within the
same line of estimates.
Relevant dependencies on parameters and special constants will be suitably emphasized using
parentheses or subscripts.\par 
As a particular case of \cite[Theorem 1.2]{Passarelli:2023:2}, the following result holds:
\begin{thm}\label{maggdiff}
    Assume that \eqref{ip1} and (\ref{ip2}) are in force for an exponent $p $ such that $ 1 < p < 2 $ and that $(u,\pi)\in W^{1,p}(\Omega,\mathbb{R}^n)\times L^{p'}(\Omega)$ is a weak solution of the system \eqref{equa} in the sense of Definition \ref{def:weak-solution}. Then, if $ f \in L^{p'}(\Omega),$ we have
  \begin{center}
     $ V_\mu(\xi) \in W^{1,2}_{loc}(\Omega, \mathbb{R}^{n \times n})$ \quad and  \quad $ \pi \in W^{\frac{2}{p'},p'}_{loc}(\Omega). $
  \end{center}
\end{thm}
\begin{flushleft}
    The next two Lemmas contain well-known technical inequalities that will be useful for our aims.
Since the value of $\mu\in[0,1]$ is fixed throughout the
article, we omit the dependence on $\mu$ in the notation.
\end{flushleft}

\begin{lem}[{\cite[Lemma 2.2]{GiaquintaModica:1986}}]\label{V-Ineq}
  For any 1 $<$ \textit{p} $<$ 2 and $\mu\in[0,1]$ we have
  \begin{equation*}
   c^{-1}(\mu^2+|\xi|^2+|\eta|^2)^{\frac{p-2}2}|\xi-\eta|^2
   \le
   |V_\mu(\xi)-V_\mu(\eta)|^2
   \le
   c(\mu^2+|\xi|^2+|\eta|^2)^{\frac{p-2}2}|\xi-\eta|^2
  \end{equation*}
  for any $\xi,\eta\in\mathbb{R}^k$ and a constant $c=c(p)>0$.
\end{lem}
\begin{lem}[{\cite[Lemma 2.1]{GiaquintaModica:1986}}]\label{Int-Ineq}
  For any 1 $<$ \textit{p} $<$ 2  and $\mu\in[0,1]$ we have
  \begin{equation*}
    c^{-1} (\mu^2+|\xi|^2+|\eta|^2)^{\frac{p-2}2}
    \le\int_0^1(\mu^2+|\xi+s(\eta-\xi)|^2)^{\frac{p-2}2}\,ds
    \le (\mu^2+|\xi|^2+|\eta|^2)^{\frac{p-2}2}
 \end{equation*}
  for any $\xi,\eta\in\mathbb{R}^k$ and a constant $c=c(p)>0$.
\end{lem}
In case $ 1 < p < 2 $ the differentiability of the function $ V_p(Du)$ implies the $W^{2,p}$ regularity of $u$. Indeed we have the following lemma.
\begin{lem}\label{differentiabilitylemma}
	Let $\Omega\subset\mathbb{R}^n$ be a bounded open set, $1<p<2$, and $v\in W^{1, p}_ {loc}\left(\Omega, \mathbb{R}^N\right)$. Then the implication
	\begin{equation*}\label{differentiabilityimplication}
		V_p\left(Dv\right)\in W^{1,2}_{loc}\left(\Omega\right) \implies v\in W^{2,p}_{loc}\left(\Omega\right) 
	\end{equation*}
	
	holds true, together with the estimate
	
	\begin{equation}\label{differentiabilityestimate}
		\int_{B_{r}}\left|D^2v(x)\right|^pdx
		\le c\cdot \left[1+\int_{B_{R}}\left|D\left(V_p\left(Dv(x)\right)\right)\right|^2+c\int_{B_R}\left|Dv(x)\right|^p\right].
	\end{equation}
	
	holds for any ball $B_R\Subset\Omega$ and $0<r<R$.
\end{lem}

The following lemma is well-known as \textit{iteration lemma},  and its
proof can be found, e.g., in \cite[Lemma 6.1, p.191]{Giusti}.
  \begin{lem}\label{lem:Giaq}
     For $R_0<R_1$, consider a bounded function
  $f:[R_0,R_1]\to[0,\infty)$ with
   \begin{equation*}
       f(r_1)\le\vartheta f(r_2)+\frac A{(r_2-r_1)^\alpha}+\frac B{(r_2-r_1)^\beta}+ C
       \qquad\mbox{for all }R_0<r_1<r_2<R_1,
    \end{equation*}
     where $A,B,C$, and $\alpha,\beta$ denote nonnegative constants
    and $\vartheta\in(0,1)$. Then we have
    \begin{equation*}
       f(R_0)\le c(\alpha,\vartheta)
    \bigg(\frac A{(R_1-R_0)^\alpha}+\frac B{(R_1-R_0)^\beta}+C\bigg).
    \end{equation*}
   \end{lem}
\subsection{Sobolev-Korn inequality}
The following result will be very  useful in the proof of our main result because it allows us, under certain assumptions, to pass from the gradient to the symmetric gradient.
\begin{lem}\label{Lemma-Korn}
  Let 1 $ < p  <  2 $ be given and
  assume that $ u\in L^p(B_\rho(x_0),\mathbb{R}^n)$ satisfies $\mathcal{E} u \in
  L^p(B_\rho(x_0),\mathbb{R}^{n \times n}_{sym})$. Then  $u\in
  W^{1,p}(B_\rho(x_0),\mathbb{R}^n)$ and
  \begin{equation}\label{korn}
    \dashint_{B_\rho (x_0)}|Du|^p\ dx\le c\dashint_{B_\rho (x_0)}|\mathcal{E} u|^p\ dx +c
    \bigg(\dashint_{B_\rho (x_0)}\Big|\frac{u-(u)_{x_0,\rho}}\rho\Big|\ dx\bigg)^p
  \end{equation}
  with a constant $c=c(n,p)$. If additionally $u=0$ on $\partial B_\rho (x_0)$, then
  \begin{equation}\label{korn-0}
    \int_{B_\rho (x_0)}|Du|^p \ dx\le c\int_{B_\rho (x_0)}|\mathcal{E} u|^p \ dx ,
  \end{equation}
  with a constant $c=c(n,p)$.
\end{lem}
A proof can be retrieved e.g. from \cite{MM}.

\subsection{A Lemma of Bogovski\u\i}
In order to let the pressure term $\nabla\pi$ in the system~\eqref{equa} disappear, it is convenient to construct
divergence free testing functions. This can be achieved by the
well-known Bogovski\u\i\ Lemma, see \cite{Bogovskii:1980}, or
\cite[Chapter 3, Section 3]{Galdi:1994}. We state it in the form needed for our purposes.
\begin{lem}\label{Bogov-lem}
  Let $B_R(x_0)$ be a ball in $\mathbb{R}^n$ and let $g\in L^p(B_R(x_0))$
  be such that $(g)_{x_0,R}=0$ and $p\in(1,\infty)$.
  Then there exists $w\in W^{1,p}_0(B_R(x_0) ,\mathbb{R}^n)$
  solving
  \begin{equation*}
   {\rm div}\, w=g\quad\mbox{in $B_R(x_0)$}
  \end{equation*}
  in the weak sense and such that
  \begin{equation*}
   \int_{B_R(x_0)}|Dw|^p\, dx\le c(n,p)\int_{B_R(x_0)}|g|^p\, dx.
  \end{equation*}
  Moreover, if the support of $g$ is contained in $B_r(x_0)$ with $0<r<R$ then also
  the support of $w$ is contained in the smaller ball $B_r(x_0)$.
\end{lem}
\subsection {Finite difference operator}

We recall some properties of the finite difference
operator that will be needed in the sequel.
 We employ the standard notation
  \begin{equation}
    \begin{gathered}
      \tau_h F(x)\equiv \tau_{h,i }F(x):= F(x+he_i)-F(x),
    \end{gathered}
\label{def:diff-quotients}
\end{equation}
for any $F\in L^1_{\textrm{loc}}(\Omega,\mathbb{R}^n)$, $i=1,\dots,n$, and
$h\neq0$ with $x+he_i,x\in\Omega$.
\par
\noindent We start with the description of some elementary properties that can be found, for example, in \cite{Giusti}.

\begin{prop}\label{findiffpr}

Let $F$ and $G$ be two functions such that
$F, G\in W^{1,p}(\Omega,\mathbb{R}^N)$, with $p\geq 1$,
and let us consider the set
$$
\Omega_{|h|}:=\left\{x\in \Omega : \mathrm{dist}(x, \partial\Omega)>|h|\right\}.
$$
Then we have 
\begin{itemize}
\item[$(i)$] $\tau_{h}F\in W^{1,p}(\Omega_{|h|},\mathbb{R}^N)$ and
$D_{i} (\tau_{h}F)=\tau_{h}(D_{i}F).$
\item[$(ii)$] If at least one of the functions $F$ or $G$ has support contained
in $\Omega_{|h|},$ then
$$
\int_{\Omega} F\, \tau_{h} G\, dx =\int_{\Omega} G\, \tau_{-h}F\, dx.
$$
\item[$(iii)$] We have
$$
\tau_{h}(F G)(x)=F(x+h)\tau_{h}G(x)+G(x)\tau_{h}F(x).
$$
\end{itemize}
\end{prop}
\noindent The next result about the finite difference operator is a kind of integral version
of the Lagrange Theorem.
\begin{lem}\label{le1} If $0<\rho<R$, $|h|<\frac{R-\rho}{2}$, $1 < p <+\infty$,
 and $F, DF\in L^{p}(B_{R})$ then
$$
\int_{B_{\rho}} |\tau_{h} F(x)|^{p}\ dx\leq c(n,p)|h|^{p} \int_{B_{R}}
|D F(x)|^{p}\ dx .
$$
Moreover, 
$$
\int_{B_{\rho}} |F(x+h )|^{p}\ dx\leq  \int_{B_{R}}
|F(x)|^{p}\ dx .
$$
\end{lem}
\begin{lem}\label{Giusti8.2}
	Let $F:\mathbb{R}^n\to\mathbb{R}^N$, $F\in L^p\left(B_R\right)$ with $1<p<+\infty$. Suppose that there exist $\rho\in(0, R)$ and $M>0$ such that
	
	$$
	\sum_{s=1}^{n}\int_{B_\rho}|\tau_{s, h}F(x)|^pdx\le M^p|h|^p
	$$
	
	for $\left|h\right|<\frac{R-\rho}{2}$. Then $F\in W^{1,p}(B_R, \mathbb{R}^N)$. Moreover
	
	$$
	\left\Vert DF \right\Vert_{L^p(B_\rho)}\le M,
	$$
	
	$$
	\left\Vert F\right\Vert_{L^{\frac{np}{n-p}}(B\rho)}\le c\left(M+\left\Vert F\right\Vert_{L^p(B_R)}\right),
	$$
	
	with $c=c(n, N, p, \rho, R)$, and 
	$$\frac{\tau_{s, h}F}{\left|h\right|}\to D_sF\qquad\mbox{ in }L^p_{loc}\left(\Omega\right),\mbox{ as }h\to0,$$ 
	for each $s=1, ..., n.$
\end{lem}

\section{Proof of Theorem \ref{Thm}}

This section is devoted to the proof of our main result that is divided in three steps. \par
In the first one, we establish an a priori estimate for the gradient of $V(\mathcal{E}u) $, assuming that such a gradient exists and belongs to $L^2_{loc}(\Omega)$. 
In the second step we use an approximation argument that allows to remove the regularity assumption on $V(\mathcal{E}u)$.  Finally we use the regularity in order to prove the regularity of the pressure.

\begin{proof}
\begin{flushleft}
 \textbf{Step 1: the a priori estimate.} Assume that $V(\mathcal{E}u) \in W^{1,2}_{loc}(\Omega,\mathbb{R}^{n \times n}_{sym})$. We fix a ball $B_R(x_0)\Subset \Omega$ and, without loss of generality, we suppose that $  0 < R < 1. $ Since the center $x_0$
  will be fixed throughout the proof, we omit it in the notation and abbreviate
  $B_R=B_R(x_0)$. 
  Since $(u,\pi)\in W^{1,p}_{loc}(\Omega,\mathbb{R}^n)\times L^{p'}_{loc}(\Omega)$ is a
  weak solution to \eqref{equa}, we have that
\begin{equation}
    \int_{\Omega} \langle a(\mathcal{E} u), \mathcal{E} \varphi \rangle \ dx =  \int_\Omega f\cdot\varphi\,dx  \quad \forall \varphi \in W_0^{1,p} : \text{div}\varphi=0.   
\end{equation}
   Since we are assumning $ f \in L^\frac{np}{n(p-1)+2-p}_{loc}(\Omega) $, it implies that  the right-hand side is well-defined. \par
We consider a cut-off function $\eta\in C_0^{\infty}(B_R) $ that will be fixed later. Given a sufficiently small increment \textit{h}, we define 
\begin{center}
    $ \textit{g} = \text{div}(\eta^2 \tau_h \textit{u}). $ \par
\end{center}
 For the construction of a divergence-free test function, we observe that
  $\mathrm{div} u=0 $ implies
\begin{equation}
    g = 2\eta \nabla \eta \cdot \tau_h u \in L^p(B_R). \notag
\end{equation}
\def\Xint#1{\mathchoice 
  {\XXint\displaystyle\textstyle{#1}}%
  {\XXint\textstyle\scriptstyle{#1}}%
  {\XXint\scriptstyle\scriptscriptstyle{#1}}%
  {\XXint\scriptscriptstyle\scriptscriptstyle{#1}}%
  \!\int}
\def\XXint#1#2#3{{\setbox0=\hbox{$#1{#2#3}{\int}$} 
  \vcenter{\hbox{$#2#3$}}\kern-.5\wd0}}
\def\dashint{\Xint-}
Moreover, we know from
  Gau\ss' theorem that
\begin{equation}
(g)_{x_0,R} = \displaystyle\dashint_{B_R} g \ dx = \dashint_{B_R} \mathrm{div}(\eta^2 \tau_h u) \ dx = \dashint_{+ \partial B_R} \eta^2 \tau_h u \cdot \textbf{n} \ dS = 0, \notag
\end{equation}
since $ \eta $ has compact support in $B_R.$
 Hence, we are in a position to apply Bogovski\u\i's
  Lemma~\ref{Bogov-lem}, which provides us with a function $\textit{w} \in W_0^{1,p}(B_R, \mathbb{R}^n) $ that has the properties:
\begin{equation}\label{choiceofw1}
    \begin{cases}

\mathrm{div}w = g \quad \mathrm{su} \ B_R, \\
\displaystyle\int_{B_R} \left| Dw \right|^p \ dx \leq c(\textit{n}, \textit{p}) \int_{B_R} \left| g \right|^p \ dx
    \end{cases}
\end{equation}
for some constant $ \textit{c} = \textit{c}(\textit{n}, \textit{p}) > 0.$ \par
\end{flushleft}
\begin{flushleft}
Let us define 
\begin{center}
    $\varphi = \tau_{-h}(\eta^2 \tau_h \textit{u}) - \tau_{-h} \textit{w} \in W_0^{1,p}(B_R, \mathbb{R}^n). $ 
\end{center}
 Since the choice of $w$ implies
  $\mathrm{div}\varphi=0$, the function $\varphi$ is admissible as test function in
 (3.1), which implies
\begin{align}
\int_{B_R} \langle a(\mathcal{E} u), \mathcal{E}(\tau_{-h}(\eta^2 \tau_h u)) \rangle \, dx 
&= \int_{B_R} \langle a(\mathcal{E} u), \mathcal{E}(\tau_{-h} w) \rangle \, dx \nonumber \\
&\quad + \int_{B_R} f \cdot \tau_{-h}(\eta^2 \tau_h u) \, dx 
- \int_{B_R} f \cdot \tau_{-h} w \, dx. \notag
\end{align}
An application of Proposition \ref{findiffpr} transforms this identity into
\begin{align}
\int_{B_R} \langle \tau_h(a(\mathcal{E} u)), \mathcal{E}(\eta^2 \tau_h u) \rangle \, dx 
&= \int_{B_R} \langle a(\mathcal{E} u), \tau_{-h}(\mathcal{E} w) \rangle \, dx \nonumber \\
&\quad + \int_{B_R} f \cdot \tau_{-h}(\eta^2 \tau_h u) \, dx 
- \int_{B_R} f \cdot \tau_{-h} w \, dx. \notag
\end{align}
Since
\begin{equation}
    \mathcal{E}(\eta^2 \tau_h u) = \eta^2 \tau_h \mathcal{E} u + 2\eta \nabla \eta \otimes \tau_h u + \eta \tau_h u \otimes \nabla \eta, \notag
\end{equation}
 we can rewrite the previous equation as
\begin{equation}
\begin{aligned}
    \int_{B_R} \eta^2 \langle \tau_h(a(\mathcal{E} u)), \tau_h \mathcal{E} u \rangle \, dx 
    &= -2 \int_{B_R} \langle \tau_h (a(\mathcal{E} u)), \eta \nabla \eta \otimes \tau_h u \rangle \, dx 
    + \int_{B_R} \langle a(\mathcal{E} u), \tau_{-h}(\mathcal{E} w) \rangle \, dx \\
    &\quad + \int_{B_R} f \cdot \tau_{-h}(\eta^2 \tau_h u) \, dx 
    - \int_{B_R} f \cdot \tau_{-h} w \, dx, \notag
\end{aligned}
\end{equation}

  where we used that $a(\mathcal{E} u)\in \mathbb{R}^{n \times n}_{sym} $ and that the symmetric
  matrices are orthogonal to the antisymmetric ones. \par
\begin{flushleft}
    With another application of Proposition \ref{findiffpr} in the first integral in the right-hand side of previous equality, we obtain that
\begin{equation}
\begin{aligned}
    \int_{B_R} \eta^2 \langle \tau_h(a(\mathcal{E} u)), \tau_h \mathcal{E} u \rangle \, dx 
    &= -2 \int_{B_R} \langle  a(\mathcal{E} u), \tau_{-h}(\eta \nabla \eta \otimes \tau_h u) \rangle \, dx \\
    &\quad + \int_{B_R} \langle a(\mathcal{E} u), \tau_{-h}(\mathcal{E} w) \rangle \, dx \\
    &\quad + \int_{B_R} f \cdot \tau_{-h}(\eta^2 \tau_h u) \, dx 
    - \int_{B_R} f \cdot \tau_{-h} w \, dx \\
    &\coloneqq I_1 + I_2 + I_3 + I_4.
\end{aligned}
\end{equation}
\end{flushleft}

Let $ \frac{R}{2}  \leq \tilde{s} < \textit{t} < \tilde{t} < R $ and fix $\eta $ so that $ \eta \in C_0^{\infty}(B_t), \eta \equiv 1 $ on $ B_{\tilde{s}}, 0 \leq \eta \leq 1,  \left| \nabla \eta \right| \leq \frac{c}{t-\tilde{s}}, $ and $\left| \nabla^2 \eta \right| \leq \frac{c}{(t-\tilde{s})^2}.$
 Thanks to  assumption (\ref{ip1}), it follows that
    \begin{equation}
        \left| I_1 \right| \leq L \int_{B_R} (\mu^2 + \left| \mathcal{E} u(x) \right|^2)^\frac{p-1}{2} \left| \tau_{-h}(\eta \nabla \eta \otimes \tau_h u) \right| \ dx.
    \end{equation}
 By (\textit{iii}) in Proposition \ref{findiffpr}, we obtain that
\begin{align}
    \tau_{-h}(\nabla \eta \cdot \eta \tau_h u) &= \tau_{-h}(\nabla \eta) \cdot \eta \tau_h u + \nabla \eta \cdot \tau_{-h}(\eta \tau_h u). \notag
    \end{align}
   and so, by the properties of $ \eta, $
   \begin{equation}
\begin{aligned}
    \left| \tau_{-h}(\nabla \eta \cdot \eta \, \tau_h u) \right| 
    &\leq \frac{c |h|}{(t - \tilde{s})^2} \, \eta |\tau_h u|
    + \frac{c}{t - \tilde{s}} \, \left| \tau_{-h}(\eta \tau_h u) \right|. 
\end{aligned}
\end{equation}
\end{flushleft}
Thus, inserting (3.5) in (3.4), it follows that
   \begin{align}
    \left| I_1 \right| 
    &\leq  \frac{c \left| h \right|}{(t - \tilde{s})^2} \int_{B_t} \eta (\mu^2 + \left| \mathcal{E} u \right|^2)^{\frac{p-1}{2}} \left| \tau_h u \right|  \ dx \\
    &\quad + \frac{c} {t - \tilde{s}} \int_{B_t} (\mu^2 + \left| \mathcal{E} u \right|^2)^{\frac{p-1}{2}} \left| \tau_{-h}(\eta  \tau_h u)) \right| \, dx. \notag
\end{align}

 H\"older's inequality with esponents $\Biggl(\textit{p}, \tfrac{p}{p-1}\Biggr)$ yields 
\begin{equation}
\begin{aligned}
\left| I_1 \right| 
&\leq \frac{c \left| h \right|}{(t-\tilde{s})^2} 
\left( \int_{B_t} (\mu^2 + \left| \mathcal{E} u \right|^2)^{\frac{p}{2}} \, dx \right)^{\frac{p-1}{p}} 
\left( \int_{B_t} \left| \tau_h u \right|^p \, dx \right)^{\frac{1}{p}} \\
&\quad + \frac{c}{t-\tilde{s}} 
\left( \int_{B_t} (\mu^2 + \left| \mathcal{E} u \right|^2)^{\frac{p}{2}} \, dx \right)^{\frac{p-1}{p}} 
\left( \int_{B_t} \left| \tau_{-h} (\eta \tau_h u) \right|^p \, dx \right)^{\frac{1}{p}} \\
&\leq \frac{c \left| h \right|^2}{(t - \tilde{s})^2} 
\left( \int_{B_t} (\mu^2 + \left| \mathcal{E} u \right|^2)^{\frac{p}{2}} \, dx \right)^{\frac{p - 1}{p}} 
\left( \int_{B_{\tilde{t}}} \left| D u \right|^p \, dx \right)^{\frac{1}{p}} \\
&\quad + \frac{c \left| h \right|}{t - \tilde{s}} 
\left( \int_{B_t} (\mu^2 + \left| \mathcal{E} u \right|^2)^{\frac{p}{2}} \, dx \right)^{\frac{p - 1}{p}} 
\left( \int_{B_{\tilde{t}}} \left| D(\eta \tau_h u) \right|^p \, dx \right)^{\frac{1}{p}},
\end{aligned}
\end{equation}
where we used Lemma \ref{le1}.

Since $\eta \tau_h \textit{u} $ has compact support, Sobolev-Korn's Inequality can be applied, thus getting 
\begin{equation}
\begin{aligned}
    \left| I_1 \right| &\leq \frac{c \left| h \right|^2}{(t - \tilde{s})^2} 
    \Biggl( \int_{B_{t}} \left( \mu^2 + \left| \mathcal{E} u \right|^2 \right)^{\frac{p}{2}} \, dx \Biggr)^{\frac{p - 1}{p}} 
    \Biggl( \int_{B_{\tilde{t}}} \left| D u \right|^p \, dx \Biggr)^{\frac{1}{p}} \\
    &\quad + \frac{c \left| h \right|}{t - \tilde{s}} 
    \Biggl( \int_{B_{t}} \left( \mu^2 + \left| \mathcal{E} u \right|^2 \right)^{\frac{p}{2}} \, dx \Biggr)^{\frac{p - 1}{p}} 
    \Biggl( \int_{B_{\tilde{t}}} \left| \mathcal{E}(\eta \tau_h u) \right|^p \, dx \Biggr)^{\frac{1}{p}}.
\end{aligned}
\end{equation}
Since 
\begin{equation}
    \mathcal{E} (\eta \tau_h u) = \eta \tau_h (\mathcal{E} u) + \nabla \eta \otimes \tau_h u, \notag
\end{equation}
arguing as in (3.5), we have
\begin{align*}
\left| I_1 \right| 
&\leq \frac{c |h|^2}{(t - \tilde{s})^2} 
\left( \int_{B_t} \left( \mu^2 + \left| \mathcal{E} u \right|^2 \right)^{\frac{p}{2}} \, dx \right)^{\frac{p - 1}{p}} 
\left( \int_{B_{\tilde{t}}} \left| D u \right|^p \, dx \right)^{\frac{1}{p}} \\
&\quad + \frac{c |h|}{t - \tilde{s}} 
\left( \int_{B_t} \left( \mu^2 + \left| \mathcal{E} u \right|^2 \right)^{\frac{p}{2}} \, dx \right)^{\frac{p - 1}{p}} 
\left( \int_{B_{t}} \left| \tau_h ( \mathcal{E} u ) \right|^p \, dx \right)^{\frac{1}{p}} \\
&\quad + \frac{c |h|}{(t - \tilde{s})^2} 
\left( \int_{B_t} \left( \mu^2 + \left| \mathcal{E} u \right|^2 \right)^{\frac{p}{2}} \, dx \right)^{\frac{p - 1}{p}} 
\left( \int_{B_{t}} \left| \tau_h u \right|^p \, dx \right)^{\frac{1}{p}}.
\end{align*}
By applying Lemma \ref{le1} to the last integral in the right-hand side of previous estimate, we obtain that
\begin{align}
\left| I_1 \right| 
&\leq \frac{c |h|^2}{(t - \tilde{s})^2} 
\left( \int_{B_t} (\mu^2 + |\mathcal{E} u|^2)^{\frac{p}{2}} \, dx \right)^{\frac{p-1}{p}} 
\left( \int_{B_{\tilde{t}}} |D u|^p \, dx \right)^{\frac{1}{p}} \nonumber \\
&\quad + \frac{c |h|}{t - \tilde{s}} 
\left( \int_{B_t} (\mu^2 + |\mathcal{E} u|^2)^{\frac{p}{2}} \, dx \right)^{\frac{p-1}{p}} 
\left( \int_{B_{t}} |\tau_h(\mathcal{E} u)|^p \, dx \right)^{\frac{1}{p}} \coloneqq I_1^1 + I_1^2.
\end{align}
By Lemma \ref{V-Ineq}, we have that:
\begin{align}
  I_1^2  
  &\leq \frac{c |h|}{t - \tilde{s}} 
  \left( \int_{B_t} \left( \mu^2 + |\mathcal{E} u|^2 \right)^{\frac{p}{2}} \, dx \right)^{\frac{p-1}{p}} \nonumber \\
  &\quad \cdot 
  \left( \int_{B_{t}} |\tau_h V(\mathcal{E} u)|^p 
  \left( \mu^2 + |\mathcal{E} u(x)|^2 + |\mathcal{E} u(x + h e_s)|^2 \right)^{\frac{p(2 - p)}{4}} \, dx \right)^{\frac{1}{p}}. \notag
\end{align}

    Since $1<p<2$,  we may use Hölder's inequality with exponents $ (\tfrac{2}{p}, \tfrac{2}{2-p}) $ to get
\begin{align}
I_1^2 
&\leq \frac{c |h|}{t - \tilde{s}} 
\left( \int_{B_t} \left( \mu^2 + |\mathcal{E} u|^2 \right)^{\frac{p}{2}} dx \right)^{\frac{p-1}{p}} 
\left( \int_{B_t} |\tau_h V(\mathcal{E} u)|^2 dx \right)^{\frac{1}{2}} \nonumber \\
&\quad \cdot 
\left( \int_{B_t} \left( \mu^2 + |\mathcal{E} u(x)|^2 + |\mathcal{E} u(x + h e_s)|^2 \right)^{\frac{p}{2}} dx \right)^{\frac{2 - p}{2p}} \nonumber \\
&\leq \frac{c |h|}{t - \tilde{s}} 
\left( \int_{B_t} \left( \mu^2 + |\mathcal{E} u|^2 \right)^{\frac{p}{2}} dx \right)^{\frac{p-1}{p}} 
\left( \int_{B_t} |\tau_h V(\mathcal{E} u)|^2 dx \right)^{\frac{1}{2}} \nonumber \\
&\quad \cdot 
\left( \int_{B_{\tilde{t}}} \left( \mu^2 + |\mathcal{E} u|^2 \right)^{\frac{p}{2}} dx \right)^{\frac{2 - p}{2p}} \nonumber \\
&\leq \frac{c |h|}{t - \tilde{s}} 
\left( \int_{B_t} |\tau_h V(\mathcal{E} u)|^2 dx \right)^{\frac{1}{2}} 
\left( \int_{B_{\tilde{t}}} \left( \mu^2 + |\mathcal{E} u|^2 \right)^{\frac{p}{2}} dx \right)^{\frac{1}{2}}, \notag
\end{align}
where we used Lemma \ref{le1}.\par
   By applying  Young's Inequality in both integrals $ I_1^1$  and $ I_1^2 $ and inserting the corresponding estimates in (3.9), it follows that
\begin{equation}
\begin{aligned}
\left| I_1 \right| 
&\leq \frac{c |h|^2}{(t-\tilde{s})^2} \int_{B_{\tilde{t}}} |D u|^p \, dx 
+ \frac{c_\sigma |h|^2}{(t - \tilde{s})^2} \int_{B_{\tilde{t}}} \left( \mu^2 + |\mathcal{E} u|^2 \right)^{\frac{p}{2}} \, dx 
+ \sigma \int_{B_{t}} \left| \tau_h V(\mathcal{E} u) \right|^2 \, dx.
\end{aligned}
\end{equation}
where the parameter $ \sigma >  0 $ will be chosen later. \par

For what concerns $I_2$, we observe that
\begin{equation}
   \left| I_2 \right| \leq L \int_{B_t} \left( \mu^2 + \left| \mathcal{E} u(x) \right|^2 \right)^{\frac{p-1}{2}} \left| \tau_{-h}(\mathcal{E} w) \right| \, dx, \notag
\end{equation}
and applying Hölder's Inequality with exponents  $(\textit{p}, \tfrac{p}{p-1}), $ we get
\begin{equation}
    \left| I_2 \right| \leq L \Biggl(\int_{B_t} (\mu^2 + \left| \mathcal{E} u \right|^2)^\frac{p}{2} \Biggr)^\frac{p-1}{p} \Biggl(\int_{B_t} \left| \tau_{-h}(\mathcal{E} w) \right|^p \ dx \Biggr)^\frac{1}{p}.\notag
\end{equation}
But since, by (3.4), it holds
\begin{equation}
    \int_{B_t} \left| \tau_{-h} (\mathcal{E} w) \right|^p \ dx \leq \int_{B_t} \left| \tau_{-h} g \right|^p \ dx, \notag
\end{equation}
we have
\begin{align}
\left| I_2 \right| 
&\leq L \left( \int_{B_t} \left( \mu^2 + \left| \mathcal{E} u \right|^2 \right)^{\frac{p}{2}} \, dx \right)^{\frac{p-1}{p}} 
\left( \int_{B_t} \left| \tau_{-h} g \right|^p \, dx \right)^{\frac{1}{p}} \notag \\
&= L \left( \int_{B_t} \left( \mu^2 + \left| \mathcal{E} u \right|^2 \right)^{\frac{p}{2}} \, dx \right)^{\frac{p-1}{p}} 
\left( \int_{B_t} \left| \tau_{-h} \left( \eta \nabla \eta \otimes \tau_h u \right) \right|^p \, dx \right)^{\frac{1}{p}}. \notag
\end{align}
Therefore the same steps taken to estimate $ I_1 $ can be applied for estimate $ I_2, $ thus getting \par
\begin{equation}
\begin{aligned}
\left| I_2 \right| 
&\leq \frac{c |h|^2}{(t-\tilde{s})^2} \int_{B_{\tilde{t}}} |D u|^p \, dx 
+ \frac{c_\sigma |h|^2}{(t - \tilde{s})^2} \int_{B_{\tilde{t}}} \left( \mu^2 + |\mathcal{E} u|^2 \right)^{\frac{p}{2}} \, dx 
+ \sigma \int_{B_{t}} \left| \tau_h V(\mathcal{E} u) \right|^2 \, dx.
\end{aligned}
\end{equation}

Let us observe that
\begin{equation}
    \left| I_3 \right| \leq \int_{B_R} \left| f \right| \left| \tau_{-h} (\eta^2 \tau_h u) \right| \ dx. \notag
\end{equation}
By the Sobolev Embedding Theorem, the a priori assumption $ V(\mathcal{E} u) \in W^{1,2}_{loc}(\Omega) $ implies $ Du \in L^{\tfrac{np}{n-2}}_{loc}(\Omega).$ Since $\eta $ has compact support, we also have that $ D(\eta u) \in L^{\frac{np}{n-2}}_{loc}(\Omega); $ as a consequence $D(\eta u) \in L^{\frac{np}{n-2+p}}_{loc}(\Omega) $ (because $ \tfrac{np}{n-2+p} < \tfrac{np}{n-2}) $. This allows to apply  Hölder’s inequality with exponents $ \Biggl(\tfrac{np}{n-2+p}, \tfrac{np}{n(p-1)+2-p}\Biggr)$ as follows
\begin{align}
\left| I_3 \right| 
&\leq \left( \int_{B_t} \left| f \right|^{\frac{np}{n(p-1) + 2 - p}} \, dx \right)^{\frac{n(p-1) + 2 - p}{np}} 
\left( \int_{B_t} \left| \tau_{-h}(\eta^2 \tau_h u) \right|^{\frac{np}{n + p - 2}} \, dx \right)^{\frac{n + p - 2}{np}} \notag \\
&\leq |h| \left( \int_{B_t} \left| f \right|^{\frac{np}{n(p-1) + 2 - p}} \, dx \right)^{\frac{n(p-1) + 2 - p}{np}} 
\left( \int_{B_{\tilde{t}}} \left| D(\eta^2 \tau_h u) \right|^{\frac{np}{n + p - 2}} \, dx \right)^{\frac{n + p - 2}{np}} \notag \\
&\leq |h| \left( \int_{B_t} \left| f \right|^{\frac{np}{n(p-1) + 2 - p}} \, dx \right)^{\frac{n(p-1) + 2 - p}{np}} 
\left( \int_{B_{\tilde{t}}} \left| \mathcal{E}(\eta^2 \tau_h u) \right|^{\frac{np}{n + p - 2}} \, dx \right)^{\frac{n + p - 2}{np}}, \notag
\end{align} 
where we used Lemma \ref{le1} and Sobolev-Korn's inequality .

Since
\begin{equation}
    \mathcal{E}(\eta^2 \tau_h u) = \eta^2 \mathcal{E}(\tau_h u) + 2\eta \nabla \eta \otimes \tau_h u = \eta^2 \tau_h(\mathcal{E} u) + 2 \eta \nabla \eta  \otimes \tau_h u, \notag
\end{equation}
from the properties of $\eta $ and an application of Lemma \ref{le1}, we obtain
\begin{align}
    \left| I_3 \right| &\leq \left| h \right| 
    \left( \int_{B_t} \left| f \right|^{\frac{np}{n(p-1)+2-p}} \, dx \right)^{\frac{n(p-1)+2-p}{np}} 
    \left( \int_{B_{t}} \left| \tau_h(\mathcal{E} u) \right|^{\frac{np}{n+p-2}} \, dx \right)^{\frac{n+p-2}{np}} \nonumber \\
    &\quad + \frac{c \left| h \right|^2}{t - \tilde{s}}  \left( \int_{B_t} \left| f \right|^{\frac{np}{n(p-1)+2-p}} \, dx \right)^{\frac{n(p-1)+2-p}{np}} 
    \left( \int_{B_{\tilde{t}}} \left| Du \right|^{\frac{np}{n+p-2}} \, dx \right)^{\frac{n+p-2}{np}}.
\end{align}

 Lemma \ref{V-Ineq} yields
\begin{align}
    \int_{B_{t}} \left| \tau_h (\mathcal{E} u) \right|^{\frac{np}{n-2+p}} \, dx  
    &\leq \int_{B_t} \left| \tau_h V(\mathcal{E}u) \right|^{\frac{np}{n-2+p}} 
    \left( \mu^2 + \left| \mathcal{E} u(x) \right|^2 + \left| \mathcal{E} u(x + h e_s) \right|^2 \right)^{\frac{np(2-p)}{4(n-2+p)}} \, dx. 
\end{align}
Since $ 1  < p  < 2, $ we are allowed to use Hölder’s inequality with the pair of conjugate exponents $\Biggl(\tfrac{2(n-2+p)}{np}, \tfrac{2(n-2+p)}{(n-2)(2-p)}\Biggr)$ in the right-hand side of previous estimate, which leads to
\begin{align}
    \int_{B_{t}} \left| \tau_h (\mathcal{E} u) \right|^{\frac{np}{n - 2 + p}} \, dx  
    &\leq \left( \int_{B_{t}} \left| \tau_h V(\mathcal{E} u) \right|^2 \, dx \right)^{\frac{np}{2(n - 2 + p)}} \nonumber \\
    &\quad \cdot \left( \int_{B_{t}} \left( \mu^2 + \left| \mathcal{E} u(x) \right|^2 + \left| \mathcal{E} u(x + h e_s) \right|^2 \right)^{\frac{np}{2(n - 2)}} \, dx \right)^{\frac{(n - 2)(2 - p)}{2(n - 2 + p)}}
\end{align}
Hence, inserting (3.14) in (3.12), we get
\begin{align}
    \left| I_3 \right| 
    &\leq |h| 
    \left( 
        \int_{B_t} \left| f \right|^{\frac{np}{n(p - 1) + 2 - p}} \, dx 
    \right)^{\frac{n(p - 1) + 2 - p}{np}} \nonumber \\
    &\quad \cdot 
    \left( 
        \int_{B_{t}} \left| \tau_h V(\mathcal{E} u) \right|^2 \, dx 
    \right)^{\frac{1}{2}} 
    \left( 
        \int_{B_{t}} \left( 
            \mu^2 + \left| \mathcal{E} u(x) \right|^2 + \left| \mathcal{E} u(x + h e_s) \right|^2 
        \right)^{\frac{np}{2(n - 2)}} \, dx 
    \right)^{\frac{(n - 2)(2 - p)}{2np}} \nonumber \\
    &\quad + \frac{c |h|^2}{t - \tilde{s}} 
    \left( 
        \int_{B_t} \left| f \right|^{\frac{np}{n(p - 1) + 2 - p}} \, dx 
    \right)^{\frac{n(p - 1) + 2 - p}{np}} 
    \left( 
        \int_{B_{\tilde{t}}} \left| D u \right|^{\frac{np}{n - 2}} \, dx 
    \right)^{\frac{n - 2}{np}}, \notag
\end{align}
where we used Holder's Inequality in the last integral with exponents $ (\tfrac{n+p-2}{n-2}, \tfrac{n+p-2}{p})  $ and Lemma \ref{le1}.
\begin{flushleft}
   We can apply Young's inequality with conjugate exponents  $ \Biggl(2, \tfrac{2p}{2-p}, \tfrac{p}{p-1}\Biggr) $ in the first term of previous estimate, thus obtaining:
\end{flushleft}
\begin{equation}
\begin{aligned}
    \left| I_3 \right| 
    &\leq c_\sigma |h|^2 \Biggl(\int_{B_R} \left| f \right|^{\frac{np}{n(p - 1) + 2 - p}} \, dx \Biggr)^\frac{n(p-1)+2-p}{n(p-1)}
    + \sigma |h|^2  \Biggl(\int_{B_{\tilde{t}}}  \left( \mu^2 + \left| \mathcal{E} u(x) \right|^2 \right)^{\frac{np}{2(n-2)}}  \, dx \Biggr)^\frac{n-2}{n} \\
    &\quad + \sigma \int_{B_t} \left| \tau_h \left( V(\mathcal{E} u) \right) \right|^2 \, dx 
    + \frac{c |h|^2}{t - \tilde{s}} 
    \left( 
        \int_{B_R} \left| f \right|^{\frac{np}{n(p - 1) + 2 - p}} \, dx 
    \right)^{\frac{n(p - 1) + 2 - p}{np}} 
    \left( 
        \int_{B_{\tilde{t}}} \left| D u \right|^{\frac{np}{n - 2}} \, dx 
    \right)^{\frac{n - 2}{np}}. \notag
\end{aligned}
\end{equation}
By  applying Young's inequality with exponents $ (\tfrac{p}{p-1}, p) $ to the last term of the right-hand side of the estimate and since $t - \tilde{s}  <  1 $ it follows that  
\begin{align*}
\left| I_3 \right| \leq\; &
\frac{c_\sigma |h|^2}{(t - \tilde{s})^{\frac{p}{p - 1}}} 
\left( \int_{B_R} \left| f \right|^{\frac{np}{n(p - 1) + 2 - p}} \, dx \right)^{\frac{n(p - 1) + 2 - p}{n(p - 1)}} 
+ \sigma |h|^2 
\left( \int_{B_{\tilde{t}}} \left( \mu^2 + \left| \mathcal{E}u(x) \right|^2 \right)^{\frac{np}{2(n - 2)}} \, dx \right)^{\frac{n - 2}{n}} \\
&+ \sigma \int_{B_t} \left| \tau_h \left( V(\mathcal{E}u) \right) \right|^2 \, dx 
+ \sigma |h|^2 \left( \int_{B_{\tilde{t}}} \left| Du \right|^{\frac{np}{n - 2}} \, dx \right)^{\frac{n - 2}{n}} \\
\leq\; &
\frac{c_\sigma |h|^2}{(t - \tilde{s})^{\frac{p}{p - 1}}} 
\left( \int_{B_R} \left| f \right|^{\frac{np}{n(p - 1) + 2 - p}} \, dx \right)^{\frac{n(p - 1) + 2 - p}{n(p - 1)}} 
+ \sigma |h|^2 
\left( \int_{B_{\tilde{t}}} \left( \mu^2 + \left| \mathcal{E}u(x) \right|^2 \right)^{\frac{np}{2(n - 2)}} \, dx \right)^{\frac{n - 2}{n}} \\
&+ \sigma \int_{B_t} \left| \tau_h \left( V(\mathcal{E}u) \right) \right|^2 \, dx 
+ c \cdot \sigma |h|^2 
\left( \int_{B_{\tilde{t}}} \left| \frac{u - (u)_{x_0, \tilde{t}}}{\tilde{t}} \right|^{\frac{np}{n - 2}} \, dx \right)^{\frac{n - 2}{n}} 
+ c \cdot \sigma |h|^2 
\left( \int_{B_{\tilde{t}}} \left| \mathcal{E}u \right|^{\frac{np}{n - 2}} \, dx \right)^{\frac{n - 2}{n}},
\end{align*}
where we applied Sobolev–Korn's inequality to the last term on the right-hand side of the  previous inequality.

Since
\begin{equation}
  \int_{B_{\tilde{t}}} \left|\mathcal{E}u \right|^\frac{np}{n-2} \ dx \leq c \int_{B_{\tilde{t}}} \left( \mu^2 + \left| \mathcal{E} u(x) \right|^2 \right)^{\frac{np}{2(n - 2)}} \, dx,  \notag
\end{equation}
then it follows that
\begin{align}
\left| I_3 \right| \leq\; &
\frac{c_\sigma |h|^2}{\left( t - \tilde{s} \right)^{\frac{p}{p-1}}} 
\left( \int_{B_R} \left| f \right|^{\frac{np}{n(p-1) + 2 - p}} \, dx \right)^{\frac{n(p-1) + 2 - p}{n(p-1)}} \nonumber \\
&\quad + c \cdot \sigma |h|^2 
\left( \int_{B_{\tilde{t}}} \left( \mu^2 + \left| \mathcal{E}u(x) \right|^2 \right)^{\frac{np}{2(n-2)}} \, dx \right)^{\frac{n-2}{n}} 
+ \sigma \int_{B_t} \left| \tau_h \left( V(\mathcal{E}u) \right) \right|^2 \, dx \nonumber \\
&\quad + c \cdot \sigma |h|^2 
\left( \int_{B_{\tilde{t}}} \left| \frac{u - (u)_{x_0, \tilde{t}}}{\tilde{t}} \right|^{\frac{np}{n-2}} \, dx \right)^{\frac{n-2}{n}} \nonumber \\
\leq\; &
\frac{c_\sigma |h|^2}{\left( t - \tilde{s} \right)^{\frac{p}{p-1}}} 
\left( \int_{B_R} \left| f \right|^{\frac{np}{n(p-1) + 2 - p}} \, dx \right)^{\frac{n(p-1) + 2 - p}{n(p-1)}} \nonumber \\
&\quad + c \cdot \sigma |h|^2 
\left( \int_{B_{\tilde{t}}} \left( \mu^2 + \left| \mathcal{E}u(x) \right|^2 \right)^{\frac{np}{2(n-2)}} \, dx \right)^{\frac{n-2}{n}} 
+ \sigma \int_{B_t} \left| \tau_h \left( V(\mathcal{E}u) \right) \right|^2 \, dx 
+ c \cdot \sigma |h|^2 \int_{B_{\tilde{t}}} \left| Du \right|^p \, dx \nonumber \\
\leq\; &
\frac{c_\sigma |h|^2}{\left( t - \tilde{s} \right)^{\frac{p}{p-1}}} 
\left( \int_{B_R} \left| f \right|^{\frac{np}{n(p-1) + 2 - p}} \, dx \right)^{\frac{n(p-1) + 2 - p}{n(p-1)}} \nonumber \\
&\quad + c \cdot \sigma |h|^2 
\left( \int_{B_{\tilde{t}}} \left| V(\mathcal{E}u) \right|^{\frac{2n}{n-2}} \, dx \right)^{\frac{n-2}{n}} 
+ \sigma \int_{B_t} \left| \tau_h \left( V(\mathcal{E}u) \right) \right|^2 \, dx 
+ c \cdot \sigma |h|^2 \int_{B_{\tilde{t}}} \left| Du \right|^p \, dx, \notag
\end{align}
where we applied Sobolev-Poincarè's inequality to the last integral of previous estimate, since $ \tfrac{np}{n-2} < \tfrac{np}{n-p}=p^* $, and we observed that 

\begin{equation}
     (\mu^2 + \left| \mathcal{E} u \right|^2)^\frac{p}{2} \leq 1+ \left| V(\mathcal{E} u) \right|^2. \notag
\end{equation}

Then, by using Sobolev's Inequality, we can conclude with this estimate
\begin{align}
\left| I_3 \right| \leq\; &
\frac{c_\sigma |h|^2}{\left( t - \tilde{s} \right)^{\frac{p}{p - 1}}} 
\left( \int_{B_R} \left| f \right|^{\frac{np}{n(p - 1) + 2 - p}} \, dx \right)^{\frac{n(p - 1) + 2 - p}{n(p - 1)}} 
+ c \cdot \sigma |h|^2 \int_{B_{\tilde{t}}} \left| D V(\mathcal{E}u) \right|^2 \, dx \nonumber \\
&\quad + c \cdot \sigma \int_{B_t} \left| \tau_h\left( V(\mathcal{E}u) \right) \right|^2 \, dx 
+ c \cdot \sigma |h|^2 \int_{B_{\tilde{t}}} \left| Du \right|^p \, dx.
\end{align}

We now observe that
\begin{equation}
    \left| I_4 \right| \leq \int_{B_t} \left| f \right| \left| \tau_{-h} w \right| \ dx = \int_{B_t} \left| f \right| \left| \tau_{-h}(\eta \nabla \eta \otimes \tau_h u) \right| \ dx. \notag
\end{equation}
Proceeding as in (3.5), we obtain the following estimate:
\begin{equation}
\begin{aligned}
\left| I_4 \right| 
&\leq \frac{c |h|}{(t - \tilde{s})^2} \int_{B_t} \left| f \right| \left| \tau_h u \right| \, dx 
+ \frac{c}{t - \tilde{s}} \int_{B_t} \left| f \right| \left| \tau_{-h} (\eta \tau_h u) \right| \, dx \\
&\leq \frac{c}{t - \tilde{s}} 
\left( \int_{B_t} \left| f \right|^{\frac{np}{n(p - 1) + 2 - p}} \, dx \right)^{\frac{n(p - 1) + 2 - p}{np}} 
\left( \int_{B_t} \left| \tau_{-h} (\eta \tau_h u) \right|^{\frac{np}{n + p - 2}} \, dx \right)^{\frac{n + p - 2}{np}} \\
&\quad + \frac{c |h|}{(t - \tilde{s})^2} 
\left( \int_{B_t} \left| f \right|^{\frac{np}{n(p - 1) + 2 - p}} \, dx \right)^{\frac{n(p - 1) + 2 - p}{np}} 
\left( \int_{B_t} \left| \tau_h u \right|^{\frac{np}{n - 2 + p}} \, dx \right)^{\frac{n - 2 + p}{np}},
\end{aligned}
\end{equation}

where we used Hölder's inequality with exponents $\Biggl(\tfrac{np}{n-2+p}, \tfrac{np}{n(p-1)+2-p}\Biggr) $. \par
By an application of Lemma \ref{le1}, we have that
\begin{equation}
\begin{aligned}
\left| I_4 \right| &\leq 
\frac{c |h|}{t - \tilde{s}} 
\left( \int_{B_t} |f|^{\frac{np}{n(p - 1) + 2 - p}} \, dx \right)^{\frac{n(p - 1) + 2 - p}{np}} 
\left( \int_{B_{\tilde{t}}} |D(\eta \tau_h u)|^{\frac{np}{n + p - 2}} \, dx \right)^{\frac{n + p - 2}{np}} \\
&\quad + \frac{c |h|^2}{(t - \tilde{s})^2} 
\left( \int_{B_t} |f|^{\frac{np}{n(p - 1) + 2 - p}} \, dx \right)^{\frac{n(p - 1) + 2 - p}{np}} 
\left( \int_{B_{\tilde{t}}} |D u|^{\frac{np}{n - 2 + p}} \, dx \right)^{\frac{n - 2 + p}{np}} \\
&\leq 
\frac{c |h|}{t - \tilde{s}} 
\left( \int_{B_t} |f|^{\frac{np}{n(p - 1) + 2 - p}} \, dx \right)^{\frac{n(p - 1) + 2 - p}{np}} 
\left( \int_{B_{\tilde{t}}} |\mathcal{E}(\eta \tau_h u)|^{\frac{np}{n + p - 2}} \, dx \right)^{\frac{n + p - 2}{np}} \\
&\quad + \frac{c |h|^2}{(t - \tilde{s})^2} 
\left( \int_{B_t} |f|^{\frac{np}{n(p - 1) + 2 - p}} \, dx \right)^{\frac{n(p - 1) + 2 - p}{np}} 
\left( \int_{B_{\tilde{t}}} |D u|^{\frac{np}{n - 2 + p}} \, dx \right)^{\frac{n - 2 + p}{np}} \\
&\leq 
\frac{c |h|}{t - \tilde{s}} 
\left( \int_{B_t} |f|^{\frac{np}{n(p - 1) + 2 - p}} \, dx \right)^{\frac{n(p - 1) + 2 - p}{np}} 
\left( \int_{B_t} |\tau_h (\mathcal{E} u)|^{\frac{np}{n + p - 2}} \, dx \right)^{\frac{n + p - 2}{np}} \\
&\quad + \frac{c |h|}{(t - \tilde{s})^2} 
\left( \int_{B_t} |f|^{\frac{np}{n(p - 1) + 2 - p}} \, dx \right)^{\frac{n(p - 1) + 2 - p}{np}} 
\left( \int_{B_t} |\tau_h u|^{\frac{np}{n + p - 2}} \, dx \right)^{\frac{n + p - 2}{np}} \\
&\quad + \frac{c |h|^2}{(t - \tilde{s})^2} 
\left( \int_{B_t} |f|^{\frac{np}{n(p - 1) + 2 - p}} \, dx \right)^{\frac{n(p - 1) + 2 - p}{np}} 
\left( \int_{B_{\tilde{t}}} |D u|^{\frac{np}{n - 2 + p}} \, dx \right)^{\frac{n - 2 + p}{np}}, \notag
\end{aligned}
\end{equation}
where in the second line of estimate we used Korn's Inequality.

Therefore, by applying Lemma \ref{le1} to the second term of the right-hand side of the previous inequality and the properties of $ \eta $, we obtain that
\begin{align}
    \left| I_4 \right| 
    &\leq \frac{c \left| h \right| }{t - \tilde{s}} 
    \left( \int_{B_t} \left| f \right|^{\frac{np}{n(p - 1) + 2 - p}} \, dx \right)^{\frac{n(p - 1) + 2 - p}{np}} 
    \left( \int_{B_{t}} \left| \tau_h(\mathcal{E}u) \right|^{\frac{np}{n + p - 2}} \, dx \right)^{\frac{n + p - 2}{np}} \nonumber \\
    &\quad + \frac{c |h|^2}{(t - \tilde{s})^2} 
    \left( \int_{B_t} \left| f \right|^{\frac{np}{n(p - 1) + 2 - p}} \, dx \right)^{\frac{n(p - 1) + 2 - p}{np}} 
    \left( \int_{B_{\tilde{t}}} \left| D u \right|^{\frac{np}{n - 2 + p}} \, dx \right)^{\frac{n - 2 + p}{np}}.\notag
\end{align}

After these manipulations, the conclusion for the estimate of $ I_4 $ is analogous to that of $ I_3. $ Therefore, by suitably applying in this specific case all the steps carried out starting from (3.13), the following estimate is obtained:
\begin{align}
\left| I_4 \right| \leq\; & 
\frac{c_\sigma |h|^2}{\left( t - \tilde{s} \right)^{\frac{2p}{p-1}}} 
\left( \int_{B_R} \left| f \right|^{\frac{np}{n(p - 1) + 2 - p}} \, dx \right)^{\frac{n(p-1)+2-p}{n(p-1)}} 
+ c \cdot \sigma |h|^2 \int_{B_{\tilde{t}}} \left| D V(\mathcal{E} u) \right|^2 \, dx \nonumber \\
& \quad + c \cdot \sigma \int_{B_t} \left| \tau_h \left( V(\mathcal{E} u) \right) \right|^2 \, dx 
+ c \cdot \sigma |h|^2 \int_{B_{\tilde{t}}} \left| Du \right|^p \, dx.
\end{align}

Inserting (3.10), (3.11), (3.15) and (3.17), applying Lemma \ref{V-Ineq} and using assumption (\ref{ip2}), we get
\begin{align}
\int_{B_{\tilde{s}}} \left| \tau_h V(\mathcal{E} u) \right|^2 \, dx 
&\leq \int_{B_t} \left( \mu^2 + \left| \mathcal{E} u(x) \right|^2 + \left| \mathcal{E} u(x + h e_s) \right|^2 \right)^{\frac{p - 2}{2}} 
\left| \tau_h \mathcal{E} u \right|^2 \eta^2 \, dx \nonumber \\
&\leq \int_{B_t} \left\langle \tau_h a(\mathcal{E} u), \tau_h (\mathcal{E} u) \right\rangle \eta^2 \, dx \nonumber \\
&\leq \left| I_1 \right| + \left| I_2 \right| + \left| I_3 \right| + \left| I_4 \right| \nonumber \\
&\leq \frac{c |h|^2}{(t - \tilde{s})^2} \int_{B_{\tilde{t}}} |D u|^p \, dx 
+ \frac{c_\sigma |h|^2}{(t - \tilde{s})^2} \int_{B_{\tilde{t}}} \left( \mu^2 + |\mathcal{E} u|^2 \right)^{\frac{p}{2}} \, dx \nonumber \\
&\quad + \frac{c_\sigma |h|^2}{\left( t - \tilde{s} \right)^{\frac{2p}{p-1}}} 
\left( \int_{B_R} \left| f \right|^{\frac{np}{n(p - 1) + 2 - p}} \, dx \right)^{\frac{n(p - 1) + 2 - p}{n(p - 1)}} \nonumber \\
&\quad + c \cdot \sigma |h|^2 \int_{B_{\tilde{t}}} \left| D V(\mathcal{E} u) \right|^2 \, dx 
+ c \cdot\sigma \int_{B_t} \left| \tau_h V(\mathcal{E} u) \right|^2 \, dx 
+ c \cdot \sigma |h|^2 \int_{B_{\tilde{t}}} \left| D u \right|^p \, dx. \notag
\end{align}
By the a priori assumption $DV(\mathcal{E}u) \in L^2_{loc}(\Omega)$ and Lemma \ref{le1}, we obtain
\begin{align}
    \int_{B_{\tilde{s}}} \left| \tau_h V(\mathcal{E} u) \right|^2 \, dx 
    &\leq \frac{\sigma |h|^2}{(t - \tilde{s})^2} \int_{B_{\tilde{t}}} |D u|^p \, dx 
    + \frac{c_\sigma |h|^2}{(t - \tilde{s})^2} \int_{B_{\tilde{t}}} \left( \mu^2 + |\mathcal{E} u|^2 \right)^{\frac{p}{2}} \, dx \nonumber \\
    &\quad + c \cdot \sigma |h|^2 \int_{B_{\tilde{t}}} \left| D V(\mathcal{E} u) \right|^2 \, dx \nonumber \\
    &\quad + \frac{c_\sigma |h|^2}{\left( t - \tilde{s} \right)^{\frac{2p}{p-1}}} 
    \left( \int_{B_R} \left| f \right|^{\frac{np}{n(p - 1) + 2 - p}} \, dx \right)^{\frac{n(p - 1) + 2 - p}{n(p - 1)}}. \notag
\end{align}
 Dividing both sides of the previous inequality by $ \left| h \right|^2, $ we obtain that
\begin{align}
\frac{1}{|h|^2} \int_{B_{\tilde{s}}} \left| \tau_h V(\mathcal{E} u) \right|^2 \, dx 
&\leq c  \cdot \sigma  \int_{B_{\tilde{t}}} \left| D V(\mathcal{E} u) \right|^2 \, dx 
+ \frac{c_\sigma}{(t - \tilde{s})^2} \int_{B_{\tilde{t}}} \left( \mu^2 + |\mathcal{E} u|^2 \right)^{\frac{p}{2}} \, dx \notag \\
&\quad + \frac{\sigma}{(t - \tilde{s})^2} \int_{B_{\tilde{t}}} |D u|^p \, dx 
+ \frac{c_\sigma}{(t - \tilde{s})^{\frac{2p}{p-1}}} 
\left( \int_{B_R} |f|^{\frac{np}{n(p - 1) + 2 - p}} \, dx \right)^{\frac{n(p - 1) + 2 - p}{n(p - 1)}}. \notag
\end{align}
Passing to the limit as $ h→0,$ by virtue of Lemma \ref{Giusti8.2} and by the a priori assumption $ V(\mathcal{E} u) \in W^{1,2}_{loc}(\Omega)$, it follows 
\begin{align}
 \int_{B_{\tilde{s}}} \left| D V(\mathcal{E} u) \right|^2 \, dx 
&\leq c  \cdot \sigma \int_{B_{\tilde{t}}} \left| D V(\mathcal{E} u) \right|^2 \, dx 
+ \frac{c_\sigma}{(t - \tilde{s})^2} \int_{B_{\tilde{t}}} \left( \mu^2 + |\mathcal{E} u|^2 \right)^{\frac{p}{2}} \, dx \notag \\
&\quad + \frac{\sigma}{(t - \tilde{s})^2} \int_{B_{\tilde{t}}} |D u|^p \, dx 
+ \frac{c_\sigma}{(t - \tilde{s})^{\frac{2p}{p-1}}} 
\left( \int_{B_R} |f|^{\frac{np}{n(p - 1) + 2 - p}} \, dx \right)^{\frac{n(p - 1) + 2 - p}{n(p - 1)}}. \notag
\end{align}
By choosing $ \sigma > 0 $ such that $ c \cdot \sigma  = \tfrac{1}{2}, $ and  $ t $  such that $ \tilde{t}  - \tilde{s} = \tfrac{1}{2}(t-\tilde{s}), $ it results that
\begin{align}
\int_{B_{\tilde{s}}} \left| D V(\mathcal{E} u) \right|^2 \, dx 
&\leq \frac{1}{2} \int_{B_{\tilde{t}}} \left| D V(\mathcal{E} u) \right|^2 \, dx 
+ \frac{c}{(\tilde{t} - \tilde{s})^2} \int_{B_{R}} \left( \mu^2 + |\mathcal{E} u|^2 \right)^{\frac{p}{2}} \, dx \nonumber \\
&\quad + \frac{c}{(\tilde{t} - \tilde{s})^2} \int_{B_{R}} |D u|^p \, dx 
+ \frac{c}{\left( \tilde{t} - \tilde{s} \right)^{\frac{2p}{p-1}}} 
\left( \int_{B_R} |f|^{\frac{np}{n(p - 1) + 2 - p}} \, dx \right)^{\frac{n(p - 1) + 2 - p}{n(p - 1)}}. \notag
\end{align}
We can apply the Iteration Lemma \ref{lem:Giaq}, thus getting
\begin{equation}\label{stimaapriori}
\begin{aligned}
\int_{B_{\frac{R}{2}}} \left| D V(\mathcal{E} u) \right|^2 \, dx 
&\leq \frac{c}{R^2} \int_{B_{R}} \left( \mu^2 + |\mathcal{E} u|^2 \right)^{\frac{p}{2}} \, dx 
+ \frac{c}{R^2} \int_{B_R} \left| D u \right|^p \, dx \\
&\quad + \frac{c}{R^{\frac{2p}{p-1}}} 
\left( \int_{B_R} \left| f \right|^{\frac{np}{n(p - 1) + 2 - p}} \, dx \right)^{\frac{n(p - 1) + 2 - p}{n(p - 1)}}.
\end{aligned}
\end{equation}
 Since by our assumption $ \textit{V}(\mathcal{E} u) \in W^{1,2}_{loc}, $ we have that $\mathcal{E} u \in W^{2,p}_{loc}, $ thus $ u \in W^{2,p}_{loc}, $ since 
\begin{equation}
    \left| D^2u \right|^p \approx \left| D(V(Du))) \right|^2 + (\mu^2 + \left| Du \right|^2)^\frac{p}{2}. \notag
\end{equation}

\textbf{Step 2: the approximation}. Now we want to remove the a priori assumption $V(\mathcal{E}u) \in W
^{1,2}_{loc}(\Omega)$, through a classical approximation argument.\par
Let us consider an open set $\Omega'\Subset\Omega$, and a function $\phi\in C^{\infty}_0(B_1(0))$ such that $0\le\phi\le1$ and $\int_{B_1(0)}\phi(x)dx=1$, and a standard family of mollifiers $\{\phi_\varepsilon\}_\varepsilon$  defined as follows
	
	\begin{equation*}
		\phi_\varepsilon(x)=\frac{1}{\varepsilon^n}\phi\left(\frac{x}{\varepsilon}\right),
	\end{equation*}
	
	for any $\varepsilon\in\left(0, d\left(\Omega', \partial\Omega\right)\right)$, so that, for each $\varepsilon$, $\phi_\varepsilon\in C^{\infty}_0\left(B_\varepsilon(0)\right)$, $0\le\phi_\varepsilon\le1$, $\int_{B_\varepsilon(0)}\phi_\varepsilon(x)dx=1.$ \par
	It is well known that, for any $h\in L^1_{loc}\left(\Omega'\right)$, setting 
	\begin{equation*}
		h_\varepsilon(x)=h\ast\phi_\varepsilon(x)=\int_{B_\varepsilon}\phi_\varepsilon(y)h(x+y)dy=\int_{B_1}\phi(\omega)h(x+\varepsilon\omega)d\omega,
	\end{equation*}
	we have $h_\varepsilon\in C^\infty\left(\Omega'\right)$.

Let us fix a ball $B_{\tilde{R}}=B_{\tilde{R}}\left(x_0\right)\Subset\Omega'$, with $\tilde{R}<1$ and, for each $\varepsilon\in\left(0, d\left(\Omega', \partial\Omega\right)\right)$, let $u_\varepsilon \in u + W^{1,p}_{0}(\Omega) $ be the unique solution to the Dirichlet problem
	\begin{equation}
    \begin{cases}
        -\mathrm{div} \left( a(\mathcal{E} u_\varepsilon) \right) + \nabla \pi_\varepsilon = f_\varepsilon & \text{in } B_{\tilde{R}}, \\
        \mathrm{div}\, u_\varepsilon = 0 & \text{in } B_{\tilde{R}}, \\
        u_\varepsilon = u & \text{on } \partial B_{\tilde{R}}, \notag
    \end{cases}
\end{equation}
	where $u\in W^{1,p}_{loc}\left(\Omega\right)$ is a weak solution to stationary $p$-Stokes system \eqref{equa}, and 
	\begin{equation}\label{fepsdef}
		f_\varepsilon=f\ast\phi_\varepsilon. \notag
	\end{equation}
	In weak form, we have
    \begin{equation}
        \int_{B_{\tilde{R}}} \langle a(\mathcal{E} u_\varepsilon), \mathcal{E} \varphi \rangle\ dx = \int_{B_{\tilde{R}}} f_\varepsilon \cdot\varphi \ dx \quad \forall \varphi \in C_0^\infty(\Omega) 
    \end{equation}
    such that $\mathrm{div} \varphi = 0.$

	Since $f\in L^{\frac{np}{n\left(p-1\right)+2-p}}_{loc}\left(\Omega\right)$, we have
	
	\begin{equation}\label{convf}
		f_\varepsilon\to f \qquad\mbox{ strongly in }L^{\frac{np}{n\left(p-1\right)+2-p}}\left(B_{\tilde{R}}\right),
	\end{equation}

	as $\varepsilon\to0.$

   For every $p > 1,$ we have that $(p^*)'=(\tfrac{np}{n-p})'=\tfrac{np}{np-n+p} < \tfrac{np}{n(p-1)+2-p}$ and so
    \begin{equation}\label{convf2}
		f_\varepsilon\to f \qquad\mbox{strongly in}L^{(p^*)'}\left(B_{\tilde{R}}\right),
	\end{equation}
	as $\varepsilon\to0.$\\

	By virtue of Theorem \ref{maggdiff}, $V(\mathcal{E} u_\varepsilon) \in W^{1,2}_{loc}({B_{\tilde{R}}})$ and so we are legitimated to apply estimates \eqref{stimaapriori}, thus getting
\begin{equation}
\begin{aligned}
\int_{B_{\frac{r}{2}}} \left| D V(\mathcal{E} u_\varepsilon) \right|^2 \, dx 
&\leq \frac{c}{r^2} \int_{B_{r}} \left( \mu^2 + |\mathcal{E} u_\varepsilon|^2 \right)^{\frac{p}{2}} \, dx 
+ \frac{c}{r^2} \int_{B_r} \left| D u_\varepsilon \right|^p \, dx \\
&\quad + \frac{c}{r^{\frac{2p}{p-1}}} 
\left( \int_{B_R} \left| f_\varepsilon \right|^{\frac{np}{n(p - 1) + 2 - p}} \, dx \right)^{\frac{n(p - 1) + 2 - p}{n(p - 1)}}.
\end{aligned}
\end{equation}
for any ball $B_r\Subset B_{\tilde{R}}$.\par

Let us observe that, since $ u_\varepsilon - u \in W_0^{1,p}(\Omega) $ and $\mathrm{div}(u_\varepsilon - u)=0 $, it is a legitimate test function in (3.19). With this choice, we have
\begin{equation}
     \int_{B_{\tilde{R}}} \langle a(\mathcal{E} u_\varepsilon), \mathcal{E}(u_\varepsilon - u) \rangle \ dx =  \int_{B_{\tilde{R}}} f_\varepsilon \cdot (u_\varepsilon - u) \ dx. \notag
\end{equation}
Therefore
\begin{equation}
    \int_{B_{\tilde{R}}} \langle a(\mathcal{E} u_\varepsilon), \mathcal{E}u_\varepsilon \rangle \ dx  = \int_{B_{\tilde{R}}} \langle a(\mathcal{E} u_\varepsilon), \mathcal{E}u \rangle \ dx + \int_{B_{\tilde{R}}} f_\varepsilon \cdot (u_\varepsilon - u) \ dx. \notag
\end{equation}
Thus, by using assumption \eqref{ip1} and the previous equality, we obtain that
\begin{align}
    \ell \int_{B_{\tilde{R}}} (\mu^2 + \left| \mathcal{E} u_\varepsilon \right|^2)^\frac{p}{2} \, dx 
    &\leq \int_{B_{\tilde{R}}} \langle a(\mathcal{E} u_\varepsilon), \mathcal{E}u_\varepsilon \rangle \, dx \notag \\
    &\leq L \int_{B_{\tilde{R}}} (\mu^2 + \left| \mathcal{E} u_\varepsilon \right|^2)^\frac{p-1}{2} \left| \mathcal{E} u \right| \, dx 
    + \int_{B_{\tilde{R}}} \left| f_\varepsilon \right| \left| u_\varepsilon - u \right| \, dx \notag \\
    &\leq \frac{\ell}{2} \int_{B_{\tilde{R}}}  (\mu^2 + \left| \mathcal{E} u_\varepsilon \right|^2)^\frac{p}{2} \, dx 
    + c \int_{B_{\tilde{R}}} \left| \mathcal{E} u \right|^p \, dx \notag \\
    &\quad + \left(\int_{B_{\tilde{R}}} \left| f_\varepsilon \right|^{(p^*)'} \, dx \right)^\frac{1}{(p^*)'} 
    \left(\int_{B_{\tilde{R}}} \left| u_\varepsilon - u \right|^{p^*} \, dx \right)^\frac{1}{p^*}, \notag
\end{align}
where we used Young and Hölder's inequalities in the last line of the previous inequality. \par
Now we use Poincarè-Sobolev's inequality in the left-hand side of previous inequality, thus getting
\begin{align}
\ell \int_{B_{\tilde{R}}} \left( \mu^2 + \left| \mathcal{E} u_\varepsilon \right|^2 \right)^{\frac{p}{2}} \, dx 
&\leq 
c \int_{B_{\tilde{R}}} \left| \mathcal{E} u \right|^p \, dx 
+ \frac{\ell}{2} \int_{B_{\tilde{R}}} \left( \mu^2 + \left| \mathcal{E} u_\varepsilon \right|^2 \right)^{\frac{p}{2}} \, dx \notag \\
&\quad + \left( \int_{B_{\tilde{R}}} \left| f_\varepsilon \right|^{(p^*)'} \, dx \right)^{\frac{1}{(p^*)'}} 
\left( \int_{B_{\tilde{R}}} \left| Du_\varepsilon - Du \right|^{p} \, dx \right)^{\frac{1}{p}}. \notag
\end{align}
Since $ u_\varepsilon - u = 0 $ on $\partial B_{\tilde{R}}$, we can use Sobolev-Korn's inequality, in order to obtain 
\begin{align}
\ell \int_{B_{\tilde{R}}} \left( \mu^2 + \left| \mathcal{E} u_\varepsilon \right|^2 \right)^{\frac{p}{2}} \, dx 
&\leq 
c \int_{B_{\tilde{R}}} \left| \mathcal{E} u \right|^p \, dx 
+ \frac{\ell}{2} \int_{B_{\tilde{R}}} \left( \mu^2 + \left| \mathcal{E} u_\varepsilon \right|^2 \right)^{\frac{p}{2}} \, dx \notag \\
&\quad 
+ \left( \int_{B_{\tilde{R}}} \left| f_\varepsilon \right|^{(p^*)'} \, dx \right)^{\frac{1}{(p^*)'}} 
\left( \int_{B_{\tilde{R}}} \left| \mathcal{E} u_\varepsilon - \mathcal{E} u \right|^p \, dx \right)^{\frac{1}{p}} \notag \\
&\leq 
c \int_{B_{\tilde{R}}} \left| \mathcal{E} u \right|^p \, dx 
+ \frac{\ell}{2} \int_{B_{\tilde{R}}} \left( \mu^2 + \left| \mathcal{E} u_\varepsilon \right|^2 \right)^{\frac{p}{2}} \, dx \notag \\
&\quad 
+ \frac{\ell}{2} \int_{B_{\tilde{R}}} \left| \mathcal{E} u_\varepsilon \right|^p \, dx 
+ \left( \int_{B_{\tilde{R}}} \left| f_\varepsilon \right|^{(p^*)'} \, dx \right)^{\frac{p'}{p^*}}, \notag
\end{align}
where we used Young's inequality. \par
Reasbsorbing the second and the third term of the right hand side of the previous inequality, we obtain
\begin{align}\label{unifbound}
    \frac{\ell}{2} \int_{B_{\tilde{R}}} (\mu^2 + \left| \mathcal{E} u_\varepsilon \right|^2)^{\frac{p}{2}} \, dx &\leq c \int_{B_{\tilde{R}}}  \left| \mathcal{E} u \right|^p \, dx 
    + \left( \int_{B_{\tilde{R}}} \left| f_\varepsilon \right|^{(p^*)'} \, dx \right)^{\frac{p'}{p^*}}.
\end{align}

Therefore, by virtue of \eqref{convf2} and  \eqref{unifbound}, the right-hand side of (3.22) can be bounded independently of $\varepsilon$. For this reason, recalling Lemma \ref{differentiabilitylemma}, we also infer that, for each $\varepsilon$, $u_\varepsilon\in W^{2, p}_{loc}\left(B_{\tilde{R}}\right)$, and recalling \eqref{differentiabilityestimate}, we also deduce that $\{u_\varepsilon\}_\varepsilon$ is bounded in $W^{2, p}_{loc}\left(B_{r}\right)$.\\
	Hence,
\begin{equation*}\label{vconvdebWp}
		u_\varepsilon\rightharpoonup v\qquad\mbox{ weakly in } W^{2,p}\left(B_{r}\right),
	\end{equation*}
	
	\begin{equation}\label{vconvforW1p}
		u_\varepsilon\to v\qquad\mbox{ strongly in } W^{1,p}\left(B_{r}\right),
	\end{equation}
	
	and
	
	\begin{equation}\label{aeconvDv}
		\mathcal{E}u_\varepsilon\to \mathcal{E}v
		\qquad\mbox{ almost everywhere in }B_r,
	\end{equation}
	
	up to a subsequence, as $\varepsilon\to0$.\\
Moreover, by the continuity of $\xi\mapsto DV(\xi)$ and \eqref{aeconvDv}, we get $DV\left(\mathcal{E}u_\varepsilon\right)\to DV\left(\mathcal{E}v\right)$ almost everywhere, and since the right-hand side of (3.22) can be bounded independently of $\varepsilon$, by Fatou's Lemma, passing to the limit as $\varepsilon\to0$ in (3.22), by \eqref{convf} and \eqref{vconvforW1p}, we get
\begin{equation}
\begin{aligned}
\int_{B_{\frac{r}{2}}} \left| D V(\mathcal{E} v) \right|^2 \, dx 
&\leq \frac{c}{r^2} \int_{B_{r}} \left( \mu^2 + |\mathcal{E} v|^2 \right)^{\frac{p}{2}} \, dx 
+ \frac{c}{r^2} \int_{B_r} \left| D v \right|^p \, dx \\
&\quad + \frac{c}{r^{\frac{2p}{p-1}}} 
\left( \int_{B_R} \left| f \right|^{\frac{np}{n(p - 1) + 2 - p}} \, dx \right)^{\frac{n(p - 1) + 2 - p}{n(p - 1)}}. \notag
\end{aligned}
\end{equation}
Our final step is to prove that $u=v$ a.e. in $B_{\tilde{R}}$.\\
	In order to prove this, we show that $v$ is a weak solution to the Dirichlet problem with boundary data $u$. It will then follow, by the uniqueness of weak solutions, that $u=v$ almost everywhere.
   First, we observe that the boundary condition is satisfied, since, by construction, $ u_\varepsilon - u \in W_0^{1,p} $, moreover $ u_\varepsilon $ converges strongly to $ v  in W^{1,p}_{loc}({B_{\tilde{R}}}) $, so also $ v-u \in W_0^{1,p} $. It follows that $ v \in  u + W_0^{1,p}({B_{\tilde{R}}}). $ \par
    Now we show that $ v $ is a weak solution of the problem. To this end, let us observe that 
\begin{align}
    & \lim_{\varepsilon \to 0} \int_{B_{\tilde{R}}} \langle a(\mathcal{E} v), \mathcal{E} \varphi \rangle \, dx \notag \\
    &= \lim_{\varepsilon \to 0} \int_{B_{\tilde{R}}} \langle a(\mathcal{E} v) - a(\mathcal{E} u_\varepsilon) + a(\mathcal{E} u_\varepsilon), \mathcal{E} \varphi \rangle \, dx \notag \\
    &= \lim_{\varepsilon \to 0} \int_{B_{\tilde{R}}} \langle a(\mathcal{E} v) - a(\mathcal{E} u_\varepsilon), \mathcal{E} \varphi \rangle \, dx + \lim_{\varepsilon \to 0} \int_{B_{\tilde{R}}} \langle a(\mathcal{E} u_\varepsilon), \mathcal{E} \varphi \rangle \, dx \notag \\
    &= \lim_{\varepsilon \to 0} \int_{B_{\tilde{R}}} \langle a(\mathcal{E} v) - a(\mathcal{E} u_\varepsilon), \mathcal{E} \varphi \rangle \, dx + \lim_{\varepsilon \to 0} \int_{B_{\tilde{R}}} f_\varepsilon \cdot \varphi \, dx \notag \\
    &= \lim_{\varepsilon \to 0} \int_{B_{\tilde{R}}} \langle a(\mathcal{E} v) - a(\mathcal{E} u_\varepsilon), \mathcal{E} \varphi \rangle \, dx \notag \\
    &\quad + \lim_{\varepsilon \to 0} \int_{B_{\tilde{R}}} (f_\varepsilon - f) \cdot \varphi \, dx + \int_{B_{\tilde{R}}} f \cdot \varphi \, dx \notag \\
    &\leq \lim_{\varepsilon \to 0} \int_{B_{\tilde{R}}} \left| \mathcal{E} v - \mathcal{E} u_\varepsilon \right|^{p-1} \left| \mathcal{E} \varphi \right| \, dx \notag \\
    &\quad + \lim_{\varepsilon \to 0} \int_{B_{\tilde{R}}} \left| f_\varepsilon - f \right|\cdot \left| \varphi \right| \, dx + \int_{B_{\tilde{R}}} f \cdot \varphi \, dx \notag \\
    &\leq \lim_{\varepsilon \to 0}  \| \mathcal{E} \varphi \|_p \left( \int_{B_{\tilde{R}}} \left| \mathcal{E} v - \mathcal{E} u_\varepsilon \right|^p \, dx \right)^{\frac{p-1}{p}} \, dx \notag \\
    &\quad + \lim_{\varepsilon \to 0} \int_{B_{\tilde{R}}} \left| f_\varepsilon - f \right| \cdot \left| \varphi \right| \, dx + \int_{B_{\tilde{R}}} f \cdot \varphi \, dx, \notag
\end{align}
where in the second-last step we used that, since $ p -2<0 $ and thanks to assumption \eqref{ip3}, 
\begin{equation}
    \left| a(\mathcal{E} v) - a(\mathcal{E}u_\varepsilon) \right| \leq c \left| \mathcal{E}v - \mathcal{E}u_\varepsilon\right|(\mu^2 + \left| \mathcal{E} v \right|^2 + \left| \mathcal{E}u_\varepsilon\right|^2)^\frac{p-2}{2} \leq c \left| \mathcal{E}v - \mathcal{E}u_\varepsilon\right| (\left| \mathcal{E}v \right|^2 + \left| \mathcal{E}u_\varepsilon \right|^2)^\frac{p-2}{2}\leq c \left| \mathcal{E}v - \mathcal{E}u_\varepsilon \right|^{p-1}. \notag
\end{equation}
    Therefore, using \eqref{vconvforW1p} and \eqref{convf}, it follows
\begin{equation}
 \int_{B_{\tilde{R}}} \langle a(\mathcal{E} v), \mathcal{E} \varphi \rangle \ dx = \int_{B_{\tilde{R}}} f \cdot \varphi \ dx \quad \forall \varphi \in C_0^\infty(\Omega) \notag 
\end{equation}
such that $\mathrm{div} \varphi = 0$. By using  the uniqueness of the weak solution, it follows that $ u = v $.
Then we can conclude that $ u $ has the regularity we were looking for and we have obtained the estimate (1.6).

\qquad \\\

\textbf{Step 3: the regularity of the pressure.} Next, we turn our attention to the differentiability of the pressure. \par Now, we choose a cut-off function $\eta \in C_0^\infty (B_{\tfrac{R}{2}}, [0,1])$ such that $\left| \nabla \eta \right| \leq \tfrac{c}{R}$; we want to prove that $\eta \pi \in W^{p-1,q'}_{loc}, $ where we denoted by $q = \tfrac{np}{n+p-2}. $ \par To this end we choose a test function $\varphi \in L^q(B_R)$ and apply Bogovski\u\i\ lemma \ref{Bogov-lem} to the function $\varphi - (\varphi)_R $. This provides us with a function $w \in W_0^{1,q}(B_R, \mathbb{R}^n) $ such that 
\begin{equation}\label{choiceoofw}
    \begin{cases}
        \mathrm{div} w= \varphi - (\varphi)_R \quad     \mathrm{in } \ B_R, \\
        \int_{B_R} \left| \mathcal{E}w \right|^q \ dx \leq c(n,p) \| \varphi \|_{L^q}^q.
    \end{cases}
\end{equation}
In order to estimate $\left| \tau_h(a(\mathcal{E}u))\right|$, we use assumption \eqref{ip3} and that $ p-2 < 0$,  with the result

\begin{equation}
    \left| a(\xi) - a(\eta) \right| \leq c \left| \xi - \eta \right|(\mu^2 + \left| \xi \right|^2 + \left| \eta \right|^2)^\frac{p-2}{2} \leq c \left| \xi - \eta \right|(\left| \xi \right|^2 + \left| \eta \right|^2)^\frac{p-2}{2} \leq c \left| \xi - \eta \right|^{p-1}. \notag
\end{equation}
for every $\xi, \eta \in \mathbb{R}^{n \times n}_{sym}.$
Then it follows that
\begin{equation}\label{tauha}
  \left| \tau_h\big(a(\mathcal{E}u)\big) \right|
  \le
  c \,  \left|\tau_h \mathcal{E}u\right|^{p-1}. 
\end{equation}

For any $h\neq0$ with $|h|\le \frac R4$, we calculate
\begin{align}\label{eta-pi-dual}
  \int_{B_{R}}\tau_h(\eta\pi)\varphi\,dx
  &=
  \int_{B_{R}}\eta\pi\tau_{-h}\varphi\,dx
  =
  \int_{B_{R}}\eta\pi\tau_{-h}\big[\varphi-(\varphi)_{R}\big]\,dx\\\nonumber
  &=
  \int_{B_{R}}\eta\pi \mathrm{div}(\tau_{-h}w)\,dx\\\nonumber
  &=
  \int_{B_{R}}\pi \mathrm{div}(\eta\tau_{-h}w)\,dx
  -
  \int_{B_{R}}\pi\nabla\eta\cdot\tau_{-h}w\,dx\\\nonumber
  &=: \mathrm{I}+\mathrm{II}.
\end{align}
Since $u$ is a weak solution of the system~\eqref{equa} in the sense
of Definition~\ref{def:weak-solution}, we can rewrite the first term to
\begin{align}
\textup{I}
&= -\int_{B_R} f \cdot \eta \, \tau_{-h}w \, dx
+ \int_{B_R} \left\langle a(\mathcal{E} u), \mathcal{E}(\eta \, \tau_{-h}w) \right\rangle dx \notag \\
&= -\int_{B_R} f \cdot \eta \, \tau_{-h}w \, dx
+ \int_{B_R} \left\langle \tau_h\big[ \eta a(\mathcal{E} u)\big], \mathcal{E} w \right\rangle dx
+ \int_{B_R} \left\langle a(\mathcal{E} u), \tau_{-h}w \otimes \nabla \eta \right\rangle dx \notag \\
&\leq \int_{B_R} \eta \, |f| \, |\tau_{-h} w| \, dx
+ \int_{B_R} \left|\tau_h\big(a(\mathcal{E} u)\big)\right| \, |\mathcal{E} w| \, dx
+ \frac{c|h|}{R} \int_{B_R} \left|a(\mathcal{E} u)\right| \, |\mathcal{E} w| \, dx \notag \\
&\quad + \frac{c}{R} \int_{B_R} \left|a(\mathcal{E} u)\right| \, |\tau_{-h} w| \, dx
\coloneqq \mathrm{I_1} + \mathrm{I_2} + \mathrm{I_3} + \mathrm{I_4}. \notag
\end{align}
where we used Proposition \ref{findiffpr} (iii), Lemma \ref{Lemma-Korn} and Lemma \ref{le1} for the last step. \par

By applying Hölder's inequality, it follows that
\begin{align}
    \left| \mathrm{I}_1 \right| 
    &\leq
    \left( \int_{B_R} \left| f \right|^{\frac{np}{n(p - 1) + 2 - p}} \, dx \right)^{\frac{n(p - 1) + 2 - p}{np}} 
    \left( \int_{B_{R}} \left| \tau_{-h}w \right|^{\frac{np}{n + p - 2}} \, dx \right)^{\frac{n + p - 2}{np}} \notag \\
    &\leq
    \left| h \right|  
    \left( \int_{B_R} \left| f \right|^{\frac{np}{n(p - 1) + 2 - p}} \, dx \right)^{\frac{n(p - 1) + 2 - p}{np}} 
    \left( \int_{B_{R}} \left| D w \right|^{\frac{np}{n + p - 2}} \, dx \right)^{\frac{n + p - 2}{np}}, \notag
\end{align}
where we used Lemma \ref{le1}.
By applying Sobolev-Korn's inequality \ref{Lemma-Korn} to the previous estimate, we obtain the following
\begin{align}\label{stimaI1}
    \left| \mathrm{I}_1 \right| 
    &\leq \left| h \right|  
    \left( \int_{B_R} \left| f \right|^{\frac{np}{n(p - 1) + 2 - p}} \, dx \right)^{\frac{n(p - 1) + 2 - p}{np}} 
    \left( \int_{B_{R}} \left| \mathcal{E} w \right|^{\frac{np}{n + p - 2}} \, dx \right)^{\frac{n + p - 2}{np}} \notag \\
    &\leq \left| h \right| \, \| f \|_{L^{q'}(B_R)} \, \| \varphi \|_{L^q(B_R)}, 
\end{align}
where we used (\ref{choiceoofw}). \par
Since $\left| \tau_h(a(\mathcal{E}u)) \right| \approx \left| \tau_h(\mathcal{E}u) \right|^{p-1}$ 
as we have already seen in (\ref{tauha}), it follows that
\begin{equation}
\left| \mathrm{I}_2 \right| 
\leq \int_{B_R} \left| \tau_h \mathcal{E} u \right|^{p-1} \left| \mathcal{E} w \right| \, dx 
\leq \left( \int_{B_R} \left| \tau_h \mathcal{E} u \right|^{p} \, dx \right)^{\frac{p-1}{p}} 
       \left( \int_{B_R} \left| \mathcal{E} w \right|^p \, dx \right)^{\frac{1}{p}}, \notag
\end{equation}
where we used  Hölder's inequality. \par
Since $u \in W^{2,p}_{loc}(\Omega), $ we can apply Lemma \ref{le1} and we obtain the following estimate
\begin{equation}\label{stimaI2}
\left| \mathrm{I}_2 \right| 
\leq \left| h \right|^{p-1} \left( \int_{B_\frac{R}{2}} \left| D(\mathcal{E} u )\right|^{p} \, dx \right)^{\frac{p-1}{p}} 
       \left( \int_{B_R} \left| \mathcal{E} w \right|^p \, dx \right)^{\frac{1}{p}} \leq  \left| h \right|^{p-1} \left( \int_{B_\frac{R}{2}} \left| D(\mathcal{E} u )\right|^{p} \, dx \right)^{\frac{p-1}{p}} 
       \| \varphi \|_{L^q(B_R)},
\end{equation}
where we used (\ref{choiceoofw}), since $ p < q. $

Let us observe that, by assumption \eqref{ip1} and since $\mathcal{E} u \in L^\frac{np}{n-2}_{loc},$ it results that $a(\mathcal{E} u) \in L^\frac{np}{(p-1)(n-2)}_{loc}.$ By applying  Hölder's inequality with exponents $ (\tfrac{np}{(n-2)(p-1)},\tfrac{np}{n+2p-2}), $ it follows that

\begin{align}\label{stimaI3}
|\mathrm{I}_3| 
&\leq \frac{c \left| h \right| }{R} 
\left( \int_{B_R} (\mu^2 + |\mathcal{E} u|^2)^{\frac{np}{2(n-2)}} \, dx \right)^{\frac{n-2}{np}} 
\left( \int_{B_R} \left| \mathcal{E} w \right|^{\frac{np}{n + 2p - 2}} \, dx \right)^{\frac{n + 2p - 2}{np}} \notag \\
&\leq \frac{c |h|}{R} 
\left( \int_{B_R} (\mu^2 + |\mathcal{E} u|^2)^{\frac{np}{2(n - 2)}} \, dx \right)^{\frac{(n - 2)(p - 1)}{np}} \| \varphi \|_{L^q(B_R)},
\end{align}
where we used (\ref{choiceoofw}), since $\tfrac{np}{n+2p-2} < q.$ \par
Arguing as in the previous step and applying Hölder's inequality, we obtain that
\begin{align}
|\mathrm{I}_4| 
&\leq \frac{c}{R} 
\left( \int_{B_R} (\mu^2 + |\mathcal{E} u|^2)^{\frac{np}{2(n-2)}} \, dx \right)^{\frac{n-2}{np}} 
\left( \int_{B_R} \left| \tau_{-h} w \right|^{\frac{np}{n + 2p - 2}} \, dx \right)^{\frac{n + 2p - 2}{np}} \notag \\
&\leq \frac{c |h|}{R} 
\left( \int_{B_R} (\mu^2 + |\mathcal{E} u|^2)^{\frac{np}{2(n - 2)}} \, dx \right)^{\frac{(n - 2)(p - 1)}{np}} 
\left( \int_{B_{R/2}} \left| D w \right|^{\frac{np}{n + 2p - 2}} \, dx \right)^{\frac{n + 2p - 2}{np}}, \notag
\end{align}

where we used Lemma \ref{le1}. \par
By appying Sobolev-Korn's inequality \ref{Lemma-Korn}, it results that
\begin{align}\label{stimaI4}
|\mathrm{I}_4| 
&\leq \frac{c |h|}{R} 
\left( \int_{B_R} (\mu^2 + |\mathcal{E} u|^2)^{\frac{np}{2(n - 2)}} \, dx \right)^{\frac{(n - 2)(p - 1)}{np}} 
\left( \int_{B_{R/2}} \left| \mathcal{E} w \right|^{\frac{np}{n + 2p - 2}} \, dx \right)^{\frac{n + 2p - 2}{np}} \notag \\
&\leq \frac{c |h|}{R} 
\left( \int_{B_R} (\mu^2 + |\mathcal{E} u|^2)^{\frac{np}{2(n - 2)}} \, dx \right)^{\frac{(n - 2)(p - 1)}{np}} 
\| \varphi \|_{L^q(B_\frac{R}{2})},
\end{align}
where we used (\ref{choiceoofw}) for the second term of the right hand side of the previous inequality, since $\tfrac{np}{n+2p-2} < q.$ \par

Since $ p-1 < 1, $ combining (\ref{stimaI1}), (\ref{stimaI2}), (\ref{stimaI3}) and (\ref{stimaI4}), we have the following estimate
\begin{align}\label{stimaI}
    \textup{I} \leq &\; |h|^{p-1} \|\varphi\|_{L^q(B_R)} \Biggl[ 
        \|f\|_{L^{q'}(B_R)} 
        + \left( \int_{B_{R/2}} |D(\mathcal{E} u)|^p \, dx \right)^{\frac{p-1}{p}} \notag \\
    &\quad + \frac{c}{R} 
        \left( \int_{B_R} (\mu^2 + |\mathcal{E} u|^2)^{\frac{np}{2(n-2)}} \, dx \right)^{\frac{(n-2)(p-1)}{np}} 
    \Biggr].
\end{align}

 Finally, we use H\"older's and Sobolev-Korn's  inequalities, Lemma \ref{le1} and (\ref{choiceoofw}) for the estimate
\begin{align*}\label{stimaII}
|\textup{II}|
&\leq \|\pi\|_{L^{q'}(B_R)} \, \frac{c}{R} 
\left( \int_{B_{R}} |\tau_{-h} w|^q \, dx \right)^{\frac{1}{q}} \\ 
&\leq \|\pi\|_{L^{q'}(B_R)} \, \frac{c |h|}{R} 
\left( \int_{B_{\frac{R}{2}}} \left| D w \right|^q \, dx \right)^{\frac{1}{q}} \\
&\leq \|\pi\|_{L^{q'}(B_R)} \, \frac{c |h|}{R} 
\| \varphi \|_{L^{q}(B_R)}.
\end{align*}

Thus, combining (\ref{stimaI}) and the previous estimate, it follows that
\begin{align}
    \left| \int_{B_R} \tau_h(\eta \pi) \, \varphi \, dx \right| 
    &\leq |h|^{p-1} \|\varphi\|_{L^q(B_R)} \Biggl[ 
        \|f\|_{L^{q'}(B_R)} 
        + \left( \int_{B_{R/2}} |D(\mathcal{E} u)|^p \, dx \right)^{\frac{p-1}{p}} \notag \\
    &\quad + \frac{c}{R} 
        \left( \int_{B_R} (\mu^2 + |\mathcal{E} u|^2)^{\frac{np}{2(n-2)}} \, dx \right)^{\frac{(n-2)(p-1)}{np}} 
        + \frac{c}{R} \| \pi \|_{L^{q'}(B_R)} 
    \Biggr]
\end{align}
for every $ \varphi \in L^q(B_R)$ and $ 0 < \left| h \right| < \tfrac{R}{4}. $ By Riesz rappresentation theorem, this implies
\begin{align}
    \Biggl( \int_{B_R} \left| \frac{\tau_h(\eta \pi)}{|h|^{p-1}} \right|^{q'} \, dx \Biggr)^{\frac{1}{q'}} 
    &\leq \|f\|_{L^{q'}(B_R)} 
    + \left( \int_{B_{R/2}} |D(\mathcal{E} u)|^p \, dx \right)^{\frac{p-1}{p}} \notag \\
    &\quad + \frac{c}{R} 
    \left( \int_{B_R} (\mu^2 + |\mathcal{E} u|^2)^{\frac{np}{2(n-2)}} \, dx \right)^{\frac{(n-2)(p-1)}{np}} 
    + \frac{c}{R} \| \pi \|_{L^{q'}(B_R)}. \notag
\end{align}
After dividing by $\left| B_R \right|^\frac{1}{q'}$, this implies that $\pi \in W^{p-1, q'}_{loc}.$

\end{proof}

\end{document}